\documentclass[11pt]{amsart}

\pdfoutput=1
\usepackage{amsmath, amssymb, amsthm, amsbsy, bbm, amsfonts, enumitem, color, verbatim, array, floatflt, yfonts, mathrsfs, xcolor}
\usepackage[margin=1.2in]{geometry}
\usepackage[all]{xy}
\usepackage[english]{babel}
\usepackage[normalem]{ulem}

\usepackage{fancyhdr} 

\usepackage{graphicx, color, overpic, xcolor}
\usepackage[pdfpagelabels]{hyperref}

\theoremstyle{plain}
\newtheorem{thm}{Theorem}[section]
\newtheorem*{thm*}{Theorem}
\newtheorem{prop}[thm]{Proposition}
\newtheorem*{prop*}{Proposition}
\newtheorem{lemma}[thm]{Lemma}
\newtheorem*{lemma*}{Lemma}
\newtheorem{corollary}[thm]{Corollary}

\theoremstyle{definition}
\newtheorem{definition}[thm]{Definition}
\newtheorem{example}[thm]{Example}

\newtheorem{notation}[thm]{Notation}

\theoremstyle{remark}
\newtheorem{remark}[thm]{Remark}

\newcommand{\R}{\mathbb{R}}
\newcommand{\Z}{\mathbb{Z}}
\newcommand{\N}{\mathbb{N}}

\newcommand{\F}{\mathbb{F}}
\newcommand{\res}{\mathbbm{k}}
\newcommand{\laurent}{\res((t))}
\newcommand{\define}{\mathrel{\mathop:}=}

\newcommand{\cO}{\mathcal O}

\newcommand{\id}{\mathbbm{1}} 
\newcommand{\App}{A} 
\newcommand{\sW}{W_0} 
\newcommand{\sS}{S_0} 
\newcommand{\aW}{W} 
\newcommand{\aS}{S} 
\newcommand{\eW}{\widetilde W} 
\newcommand{\Cf}{\mathcal{{C}}_{f}} 
\newcommand{\fa}{{\bf{c}_f}} 
\newcommand{\Cfm}{\mathcal{{C}}^{op}_f} 
\newcommand{\Cw}{\mathcal{C}} 
\newcommand{\aH}{\mathcal{H}} 

\newcommand{\type}{\mathrm{type}} 

\newcommand{\Conv}{\operatorname{Conv}} 

\newcommand{\Sh}{\operatorname{Sh}} 
\renewcommand{\star}{\operatorname{star}} 
\newcommand{\CAT}{\operatorname{CAT}}

\newcommand{\x}{\mathbf{x}}
\newcommand{\y}{\mathbf{y}}
\newcommand{\ys}{\mathbf{ys}}
\newcommand{\z}{\mathbf{z}}
\newcommand{\w}{\mathbf{w}}

\renewcommand{\L}{L}

\numberwithin{equation}{subsection}

\definecolor{amethyst}{rgb}{0.6, 0.4, 0.8}

\begin{document}

\hypersetup{pdfauthor={Milicevic, Naqvi, Schwer, Thomas},pdftitle={Chimney Retractions}}

\title[A gallery model for affine flag varieties]{A gallery model for affine flag varieties\\ via chimney retractions}

\author{Elizabeth Mili\'{c}evi\'{c}, Yusra Naqvi, Petra Schwer, and Anne Thomas}
\address{Elizabeth Mili\'{c}evi\'{c}, Department of Mathematics \& Statistics, Haverford College, 370 Lancaster Avenue, Haverford, PA, USA
\newline Yusra Naqvi, School of Mathematics \& Statistics, Carslaw Building F07,  University of Sydney NSW 2006, Australia
\newline Petra Schwer, Department of Mathematics, Universitätsplatz2, Otto-von-Guericke University of Magdeburg, Germany
\newline Anne Thomas, School of Mathematics \& Statistics, Carslaw Building F07,  University of Sydney NSW 2006, Australia}
\email{emilicevic@haverford.edu, yusra.naqvi@sydney.edu.au, petra.schwer@ovgu.de \newline anne.thomas@sydney.edu.au}

\thanks{The first author was partially supported by NSF Grant DMS \#1600982.}
\thanks{The second author and this research was supported in part by ARC Grant DP180102437.}
\thanks{The third author was partially supported by the DFG Project SCHW 1550/4-1.}
\date{ \today }


\begin{abstract}
This paper provides a unified combinatorial framework to study orbits in certain affine flag varieties via the associated Bruhat--Tits buildings.  We first formulate, for arbitrary affine buildings, the notion of a chimney retraction.  This simultaneously generalizes the two well-known notions of retractions in affine buildings: retractions from chambers at infinity and retractions from alcoves.  We then present a recursive formula for computing the images of certain minimal galleries in the building under chimney retractions, using purely combinatorial tools associated to the underlying affine Weyl group.  Finally, for Bruhat--Tits buildings in the function field case, we relate these retractions and their effect on minimal galleries to double coset intersections in the corresponding affine flag variety. 
\end{abstract}

\maketitle

\section{Introduction}\label{sec:intro}

Buildings enable the study of the groups which act on them via a combination of algebraic, combinatorial, and geometric methods.  The current work exemplifies this combination of approaches, in the setting of affine buildings.  An affine building $X$ of type $(\aW,\aS)$ is a union of subcomplexes called apartments.  Each apartment is a copy of the Coxeter complex for $(\aW,\aS)$; that is, the tesselation of Euclidean space whose maximal simplices, called alcoves, are in bijection with the elements of $\aW$.  Simplicial maps, called retractions, from $X$ to a fixed apartment allow us to study structures in the entire building by considering (their image in) just one apartment.  

We begin by simultaneously generalizing the two well-known classical notions of retractions of an affine building, by formalizing retractions from chimneys; see Section~\ref{sec:chimneysRetractionsGalleries}.  We then establish a combinatorial recursion which allows us to compute the images of certain galleries in $X$ under chimney retractions; a special case of this recursion is stated as Theorem~\ref{thm:recursionIntro} later in this introduction.  Many (but not all) affine buildings are of algebraic origin, in the sense that they are the Bruhat--Tits building for a reductive group $G$ over a local field~$F$ with a discrete valuation; see \cite{BruhatTits2, Weiss-affine}.   In this setting, a  chimney corresponds to a choice of a (spherical) standard parabolic subgroup of $G(F)$ and an element of the affine Weyl group~$\aW$. In Theorems~\ref{thm:DoubleCosetIntersectionI} and~\ref{thm:Grassmannian-intro}, for $F$ a function field we relate intersections of certain double cosets in the affine flag variety and affine Grassmannian, respectively, to families of galleries obtained by applying the corresponding chimney retraction.  The case of general partial flag varieties is treated in Section~\ref{sec:galleriesCosets}; see Theorem~\ref{thm:DoubleCosetIntersectionParahoric}.  The later sections of this introduction are mostly dedicated to formally stating these results and giving some examples.  We provide many more examples in Section~\ref{sec:examples}.

\begin{figure}[ht]
\begin{center}
\resizebox{0.5\textwidth}{!}
{
\begin{overpic}{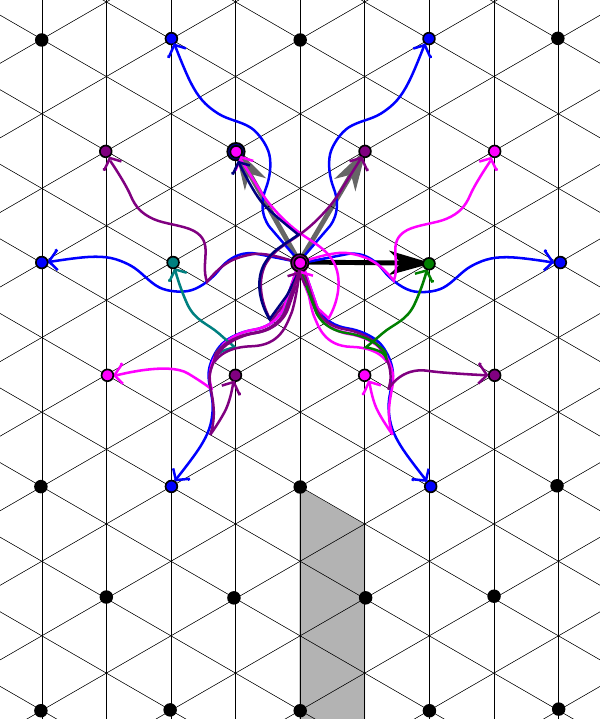}
\put(52,65){$\alpha$}
\put(30,73){$\beta$}
\put(60,96){$\lambda$}
\end{overpic}
}
\caption{The collection of all galleries of type $\vec{\lambda}$ which are positively folded with respect to the chimney represented by the gray shaded region.}
\label{fig:chimneyIntro}
\end{center}
\end{figure}

Our approach to the proof of Theorems~\ref{thm:DoubleCosetIntersectionI} and~\ref{thm:Grassmannian-intro} 
 is a generalization of that in~\cite{PRS}, which relies on the root group structure of $G(F)$, and the identification of the (points of the) varieties in question with certain simplicial substructures of the building.  The various double cosets then become orbits of simplices which can be described as pre-images under certain retractions, or equivalently as certain sets of positively folded galleries.  We note that this approach is independent of the characteristic of $F$.  This fact suggests that one should think of the Bruhat--Tits building as a tool to simultaneously study affine Grassmannians and other (partial) affine flag varieties in a characteristic-free setting.

Part of the motivation for this work is the following application to affine Deligne--Lusztig varieties.  Consider the case where $F = \res((t))$ is the field of Laurent series with $\res$ an algebraic closure of a finite field, and let $G$ be a split connected reductive group over $\res$.  In this setting, a special case of Theorem~\ref{thm:DoubleCosetIntersectionI} was combined with Theorem~6.3.1  of~\cite{GHKR} to relate positively folded galleries in the Bruhat--Tits building for $G(F)$ to affine Deligne--Lusztig varieties in the affine flag variety parameterized by translations, in a manner which makes the corresponding nonemptiness and dimension calculations tractable; see Theorem~5.8 of~\cite{MST1}.  
An analogous statement relating more general affine Deligne--Lusztig varieties to positively folded galleries can likewise be obtained by combining Theorem~\ref{thm:DoubleCosetIntersectionI} with Theorem~11.3.1 of~\cite{GHKRadlvs}.

\subsection{Retractions and gallery combinatorics}

Positively folded galleries were first introduced by Littelmann and  Gaussent--Littelmann in \cite{LittelmannLR, GaussentLittelmann} in the context of understanding highest weight representations of complex semisimple algebraic groups. Ram then developed the machinery of positively folded alcove walks~\cite{Ram}, which provides a slightly different approach to the same underlying combinatorics.  In this paper, we develop the more general combinatorics of galleries which are positively folded with respect to a chimney in an arbitrary affine building $X$; see Sections~\ref{sec:preliminaries} and~\ref{sec:chimneysRetractionsGalleries} for  notation and formal definitions.  All of these positively folded galleries are the images in a single apartment $\App$ of certain galleries in the entire building $X$, under the application of certain retractions $X \to \App$.

A \emph{chimney} in an apartment $A$ is a point in the boundary of $\App$ which is represented by a region of $\App$ lying between pairs of adjacent $\alpha$-hyperplanes for all positive roots $\alpha$ in some sub-root system, and which ``goes off to infinity" in the other root directions.  In the special case that the sub-root system equals the entire root system, an associated chimney is represented by an alcove in $\App$, rather than a boundary point.  An example of a region representing a chimney where $(W,S)$ is of type $\tilde A_2$ is the gray shaded strip in Figure~\ref{fig:chimneyIntro}.  The chimney here has associated sub-root system $\{ \pm \alpha\}$, with $\alpha$ the black (positive) root, and the other positive roots shown in gray.  Chimneys were first introduced by Rousseau in~\cite{Rousseau77} (although our definition is different), and they are related to the generalized sectors for affine buildings appearing in~\cite{CapraceLecureux,Charignon}, for instance.

A \emph{retraction from a chimney} folds the entire building $X$ away from the chimney down onto~$\App$.  Retractions from chimneys were, so far as we know, first described in~\cite{GHKRadlvs} (for certain Bruhat--Tits buildings), and the retractions from~\cite{GHKRadlvs} were used in~\cite{HKM}.  Neither the existence of chimney retractions nor our methods of proof will surprise experts.  Chimney retractions specialize to both of the classical retractions of affine buildings onto an apartment, namely the retraction centered at an alcove inside the apartment, and the retraction centered at a chamber at infinity; see for example \cite{AB}. 

A gallery in $\App$ is \emph{positively folded with respect to the chimney} if all its folds occur on sides of hyperplanes facing away from the chimney. The corresponding orientation on hyperplanes induced by a chimney generalizes the periodic orientations considered in~\cite{PRS, MST1} (see also~\cite{GraeberSchwer} for a general treatment of orientations).  The galleries in Figure~\ref{fig:chimneyIntro} are the smooth paths starting at the origin and ending in arrows, and all of these galleries are positively folded with respect to the chimney represented by the gray shaded region.

\subsection{Recursive description of shadows}
 A \emph{shadow} with respect to a chimney is the set of end-simplices of all positively folded galleries of a particular (fixed) type.  In Figure~\ref{fig:chimneyIntro}, the shadow of the coroot lattice element $\lambda$ is the set of all end-vertices of the depicted galleries; that is, the set of (colored) dots in which an arrow ends.  Shadows were first defined in~\cite{GraeberSchwer}.   Our first main result, Theorem~\ref{thm:AlcoveRecursion}, establishes a recursive formula for shadows with respect to general chimneys in the setting of an arbitrary affine building of type $(W,S)$; a special case of Theorem~\ref{thm:AlcoveRecursion} is stated as Theorem~\ref{thm:recursionIntro} below.  

Let $y \in \aW$ and write $\y$ for the corresponding alcove in a fixed apartment $\App$ of $X$.  We denote by  $\Sh_{\y}(\lambda)$ the shadow of $\lambda$ with respect to the chimney for~$\y$; i.e.~the set of end-vertices of galleries in $\App$ of a fixed type $\vec{\lambda}$ which are folded away from $\y$.  For $s \in S$, let $H$ be the hyperplane separating the alcoves $\y$ and $\ys$, and let $r_H$ be the reflection of $\App$ in $H$.  This setup permits the following recursive description of the shadow $\Sh_{\ys}(\lambda)$.

\begin{thm}\label{thm:recursionIntro}  Let $\lambda$ be any element of the coroot lattice.  For any $y \in W$ and $s \in S$ such that $\ell(ys) > \ell(y)$, we have
\[
\Sh_{\ys}(\lambda) = \Sh_\y(\lambda) \cup r_H(H^{ys} \cap \Sh_\y(\lambda)),
\]
 where $H^{ys}$ is the half-apartment which is bounded by $H$ and contains $\ys$.
\end{thm}

\begin{figure}[ht]
\begin{center}
\resizebox{0.5\textwidth}{!}
{
\begin{overpic}{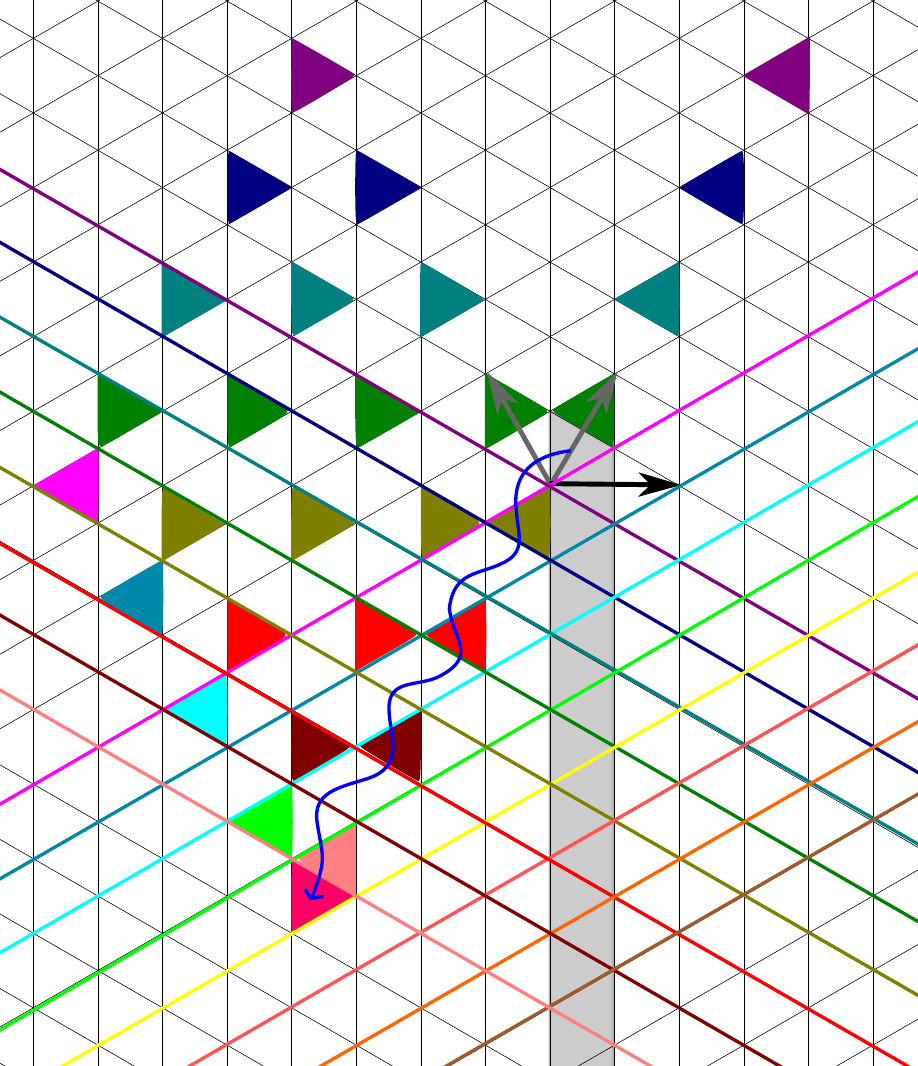}
\put(58.5,55.5){$\alpha$}
\put(28,14){$\x$}
\put(52,59){$\fa$}
\end{overpic}
}
\caption{The colored-in alcoves are the shadow of the alcove $\x$ (in pink) with respect to the chimney represented by the gray shaded region.}
\label{fig:AlcoveShadowA2}
\end{center}
\end{figure}

The idea of the recursion is that $r_H$ folds across the hyperplane $H$ all parts of those galleries from the previous step which protrude across $H$, and thus folds these galleries away from the alcove $\ys$.  For example, we used Theorem~\ref{thm:AlcoveRecursion} repeatedly to compute the shadow appearing in Figure~\ref{fig:AlcoveShadowA2}.  In this figure, a hyperplane $H$ which is crossed as one moves from the base alcove $\fa$ along the shaded region representing the chimney is shown in the same color as those alcoves in the shadow which are (first) obtained by applying the reflection $r_H$.  See Theorem~\ref{thm:AlcoveRecursion} for the general version of Theorem \ref{thm:recursionIntro}, and its proof.  Figure~\ref{fig:chimneyIntro} shows that the shadow of~$\lambda$ is a (proper) subset of the intersection of the coroot lattice with the $\lambda$-Weyl polytope, a classical result which we recover from our recursion in Corollary~\ref{cor:WeylPolytope}.

Recursions similar to Theorem \ref{thm:recursionIntro} have appeared in~\cite{GoertzHeDim, Be2}, by adapting an algebraic reduction method of Deligne and Lusztig~\cite{DL} to affine flag varieties.  Different algorithms for computing certain shadows using word combinatorics can be found in~\cite{GraeberSchwer}, and those algorithms could be generalized to the cases considered in this work.  

As we explain further after the proof of Theorem~\ref{thm:AlcoveRecursion}, determining the multiplicities of elements of the shadow is a delicate question. This can already be seen from Figure~\ref{fig:chimneyIntro}, where the origin and the element $\beta = \beta^\vee$ of the coroot lattice have multiplicity $2$, since there are two distinct galleries of type $\vec{\lambda}$ ending at these points, but all other vertices in the shadow of~$\lambda$ have multiplicity $1$.  We leave the determination of such multiplicities, and their algebraic interpretation, to future work.

\subsection{Algebraic interpretation of shadows}

Shadows of chimneys also relate sets of positively folded galleries to certain orbits in (partial) affine flag varieties over function fields.  The special cases of the affine flag variety and the affine Grassmannian are highlighted below in Theorems~\ref{thm:DoubleCosetIntersectionI} and~\ref{thm:Grassmannian-intro}, respectively, which represent the next main results in this paper. It is well-known that the nonemptiness patterns of such double coset intersections are quite delicate, at least in the case of the affine flag variety.  Recursive formulas such as Theorem~\ref{thm:recursionIntro} can be used, however, to precisely determine the sets of folded galleries arising in Theorems~\ref{thm:DoubleCosetIntersectionI} and~\ref{thm:Grassmannian-intro}.   

We now let $F = \laurent$, for any field $\res$, and let $G$ be any split, connected, reductive group over $F$.  Let $\cO = \res[[t]]$ be the ring of integers of $F$.  Let $X$ be the Bruhat--Tits building for $G(F)$, with standard apartment $\App$.  Write $I$ for the Iwahori subgroup of $G(F)$.  Note that $I$ is a subgroup of $G(\cO)$, and that $I$ is the stabilizer in $G(F)$ of the base alcove of $\App$.  Put $K = G(\cO)$, and note that $K$ is the stabilizer in $G(F)$ of the origin of $\App$. Denote by $P = P(F)$ any standard spherical parabolic subgroup of $G(F)$, with Levi decomposition $P = MN = M(F)N(F)$, where $M$ is the Levi component and $N$ is its unipotent radical. Following~\cite{GHKRadlvs}, define the subgroup $I_P$ of $G(F)$ by  $$I_P = (I \cap M)N = (I \cap M(F))N(F).$$ 
Given an element of the affine Weyl group $y \in \aW$, denote by $(I_P)^y$ the conjugate $y I_P y^{-1}$.  

The \emph{$P$-chimney} lies between the hyperplanes $H_{\alpha}$ and $H_{\alpha,1}$ of $\App$ for all positive roots $\alpha$ in~$M$,  and is antidominant for all roots in $N$. The \emph{$(P,y)$-chimney} is then the image of the $P$-chimney under the action of $y$.  For example, in Figures~\ref{fig:chimneyIntro} and~\ref{fig:AlcoveShadowA2}, if $M$ has root system~$\{ \pm \alpha \},$ then the shaded region represents the $P$-chimney.  A \emph{labeled folded alcove walk} is a folded gallery in which certain simplices have been labeled by elements of the residue field $\res$ of $F$; see Section~\ref{sec:galleriesCosets} for precise definitions.

\begin{thm}\label{thm:DoubleCosetIntersectionI}  For any $x, y, z \in \aW$, there is a bijection between the points of the intersection 
	\[  I x I  \cap (I_P)^y z I \] 
	and the set of labeled folded alcove walks of type $\vec{x}$ from the base alcove to $\z$ which are positively folded with respect to the $(P,y)$-chimney. 
\end{thm}

In the next statement, we write $\sW$ for the spherical Weyl group, $t^\lambda$ for the translation by the element~$\lambda$ of the coroot lattice, and $x_\lambda \in \aW$ for the minimal length representative of the coset $t^\lambda \sW$ in $\aW/\sW$.  Note that any double coset $K t^\lambda K$ has a dominant representative. 

\begin{thm}\label{thm:Grassmannian-intro}  Let $\lambda$ and $\mu$ be in the coroot lattice, with $\lambda$ dominant, and let $y \in \aW$.  Then there is a bijection between the points of the intersection 
	\[Kt^\lambda K \cap (I_P)^y t^\mu K\] 
and the union over $w \in \sW$ of the set of labeled folded alcove walks of type $\overrightarrow{wx_\lambda}$ from the base alcove to an alcove containing $\mu$ which are positively folded with respect to the $(P,y)$-chimney. 
\end{thm}

\noindent In Section~\ref{sec:galleriesCosets}, we prove Theorem~\ref{thm:DoubleCosetIntersectionI} as Theorem~\ref{thm:doublecosety}, and then use this to prove a generalization of Theorems~\ref{thm:DoubleCosetIntersectionI} and~\ref{thm:Grassmannian-intro} in which $I$ and $K$, respectively, are replaced by parahoric subgroups of $G(F)$ which stabilize faces of the base alcove containing the origin; see Theorem \ref{thm:DoubleCosetIntersectionParahoric}.  We also establish a general nonemptiness result for double coset intersections in Theorem~\ref{thm:IntersectionParahoricShadow}.

In some special cases, (parts of the statements of) Theorems~\ref{thm:DoubleCosetIntersectionI} and~\ref{thm:Grassmannian-intro} are folklore or recover earlier results.  For instance, if $P = B$, the Borel subgroup of $G(F)$, then  Theorem~7.1 of~\cite{PRS} is the corresponding case of Theorem~\ref{thm:DoubleCosetIntersectionI}. More generally, for any $w$ in the spherical Weyl group, Corollary~5.5 of~\cite{MST1} is the case of Theorem~\ref{thm:DoubleCosetIntersectionI} where the intersection under consideration is $IxI \cap U^w z I$.  Inspired by Peterson's work on affine Schubert calculus \cite{Pet}, in joint work with Ram, the first author proved a bijection similar to Theorem~\ref{thm:DoubleCosetIntersectionI} on intersections of positive and negative Iwahori orbits in the affine flag variety~\cite{MR}. Special cases of the orbit intersections appearing in Theorem~\ref{thm:Grassmannian-intro} were studied in \cite[Thm 2.6.11(3)-(4)]{MacdonaldSpherical}, \cite[Paragraph 4.4]{BruhatTits}, and~\cite{GaussentLittelmann}.  
In~\cite{Hitzel}, the third author used similar double cosets to prove an analog for Bruhat--Tits buildings of Kostant's convexity theorem \cite{Kostant} for symmetric spaces, which itself is a variant of a classical result of Schur \cite{Schur} and Horn \cite{Horn}.  
A count similar to the statement of Theorem~\ref{thm:DoubleCosetIntersectionParahoric} appears in the work of Abramenko, Parkinson, and Van Maldeghem \cite{APvM}, where intersections of Weyl-distance spheres in arbitrary locally finite buildings are studied (see Remark~\ref{rem:APVM} for further details).

Figures~\ref{fig:chimneyIntro} and~\ref{fig:AlcoveShadowA2} illustrate the statements of Theorems~\ref{thm:Grassmannian-intro}  and~\ref{thm:DoubleCosetIntersectionI}, respectively, for $y=1$.  If $\mu$ is the end-vertex of a gallery in Figure~\ref{fig:chimneyIntro}, equivalently $\mu$ is in the shadow of~$\lambda$ with respect to the $P$-chimney, then by Theorem~\ref{thm:Grassmannian-intro}  any labeling of a gallery in this figure which ends at $\mu$ corresponds bijectively with a point in the (nonempty) intersection $K t^\lambda K \cap I_P t^\mu K$, and conversely.   In Figure~\ref{fig:AlcoveShadowA2}, any labeling of a gallery with final alcove $\z$ in this shadow corresponds to a point in the (nonempty) intersection $IxI \cap I_P z I$, and conversely.

\subsection{Examples and code}

We conclude this work with Section~\ref{sec:examples}, which provides many examples of shadows in rank~$2$.  We observe some features of these shadows and illustrate how to use the recursion from Theorem~\ref{thm:AlcoveRecursion}.  
We have written Maple code for computing shadows in types $\tilde{A}_2$, $\tilde{C}_2$, and $\tilde{G}_2$, and this code is available from the authors upon request.

\subsection*{Acknowledgements}  We thank Ulrich G\"ortz, Tom Haines, Jean L\'ecureux, Peter Littelmann, James Parkinson, Arun Ram, and Geordie Williamson for helpful conversations and correspondence.  We also thank the University of Sydney Mathematical Research Institute for hosting a visit by the third author. The authors are grateful to the anonymous referees for their careful reading which resulted in several key clarifications.

\section{Preliminaries}\label{sec:preliminaries}

We now recall the definitions we will need and fix notation.  We assume knowledge of affine Coxeter systems and affine buildings at the level of the references~\cite{AB},~\cite{Humphreys}, or~\cite{Ronan}.  

In order to avoid notational complexities, we will assume for the remainder of this work that the affine Coxeter system $(W,S)$ is irreducible.  Our results in Sections~\ref{sec:chimneysRetractionsGalleries} and~\ref{sec:recursion} can be extended to reducible affine Coxeter systems $(W,S)$, for example where $(W,S)$ is of type $\tilde A_1 \times \tilde A_1$ so that the associated building $X$ can be realized as a product of trees.  We leave this extension as an exercise for the reader.

\subsection{Affine Coxeter systems and affine buildings}\label{sec:buildings} 

Our approach in this section is very similar to that of~\cite{LMPS}, which provides a way to discuss the affine Coxeter system $(\aW,\aS)$ without assuming that $(\aW,\aS)$ is the affine Weyl group for any particular group $G(F)$.  We do this so that we can then discuss an arbitrary affine building $X$ of type $(\aW,\aS)$.  

Let $(\aW,\aS)$ be an irreducible affine Coxeter system of rank $n+1$ and let $V$ be the associated $n$-dimensional real vector space on which $\aW$ acts, which we can identify with $n$-dimensional Euclidean space.  Write $v_0$ for the origin of $V$.  Then $(\aW,\aS)$ has associated spherical Coxeter system $(\sW,\sS)$ such that $\sW$ is the stabilizer in $\aW$ of $v_0$ and $\sS = \{ s_1, \dots, s_n \}$ is the set of elements of $\aS = \{ s_0,s_1, \dots, s_n\}$ which fix $v_0$.  We typically use the letters $x$, $y$, and $z$ for elements of $\aW$ and $u$, $v$, and $w$ for elements of $\sW$.  For algebraic reasons to be found below, we denote the elements of the translation subgroup of $\aW$ by $\{ t^\lambda \mid \lambda \in R^\vee\}$, where $R^\vee$ is a certain lattice in $V$ (which is preserved by the actions of $\aW$ and $\sW$).  For all $\lambda, \mu \in R^\vee$, we have $t^\lambda t^\mu = t^{\lambda + \mu} = t^{\mu + \lambda} = t^\mu t^\lambda$.  Any element $x$ of $\aW$ can be expressed uniquely as $x = t^\lambda w$ where $\lambda \in R^\vee$ and $w \in \sW$, and given any $u \in \sW$ and $\lambda \in R^\vee$, we have $u t^\lambda u^{-1} = t^{u\lambda}$.  We write $\ell:\aW \to \Z$ for the length function of the Coxeter system $(\aW,\aS)$.

Let $\Phi$ be a crystallographic root system such that $(\sW,\sS)$ is the Weyl group of $\Phi$ and $R^\vee \subset V$ is the associated coroot lattice.  Denote the positive roots in $\Phi$ by $\Phi^+$, and denote the positive simple roots in $\Phi$ by $\Delta = \{ \alpha_1,\dots,\alpha_n\}$.  Given $J \subseteq \{1,\dots,n\}$, we denote by $W_J$ the standard (spherical) parabolic subgroup of $\sW$ generated by $\{ s_j \mid j \in J \}$.  Recall that $(W_J, \{ s_j\}_{j\in J})$ is a Coxeter system.  We will write $\Phi_J$ for the sub-root system of $\Phi$ associated to $W_J$, with positive roots $\Phi_J^+ = \Phi^+ \cap \Phi_J$ and basis of positive simple roots $\Delta_J = \{ \alpha_j \mid j \in J \}$.  In particular, if $J = \emptyset$ then $W_J$ is trivial and $\Phi_J = \Delta_J= \emptyset$, while if $J = \{ 1,\dots,n\}$ then $W_J = \sW$, $\Phi_J = \Phi$, and $\Delta_J = \Delta$.

For each $\alpha \in \Phi$ and $k \in \Z$, we write $H_{\alpha,k}$ for the affine hyperplane or \emph{wall} of $V$ given by $H_{\alpha,k} = \{ v \in V \mid \langle \alpha, v \rangle = k \}$ where the brackets denote the standard inner product on $V$. Note that $H_{\alpha,k} = H_{-\alpha,-k}$ for all $\alpha \in \Phi$ and $k \in \Z$.  Write $s_{\alpha,k}$ for the (affine) reflection fixing $H_{\alpha,k}$ pointwise.  Then each $s_{\alpha,k}$ is an element of $\aW$, and every reflection in $\aW$ is of this form.  We will sometimes put $H_{\alpha} = H_{\alpha,0}$ and $s_{\alpha} = s_{\alpha,0}$.  We remark that $s_i = s_{\alpha_i}$ for $1 \leq i \leq n$, and $s_0 = s_{\widetilde \alpha,1}$ where $\widetilde\alpha$ is the (unique) highest root in $\Phi$.   
Let $\aH = \{ H_{\alpha,k} \mid \alpha \in \Phi, k \in \Z\}$. The group $\aW$ acts on the set $\aH$ as follows.  For $\alpha \in \Phi$, $k \in \Z$, $\mu \in R^\vee$, and $u \in \sW$, we have $t^\mu u \cdot H_{\alpha,k} = H_{u\alpha,k + \langle u\alpha, \mu \rangle}$. 

An \emph{alcove} is the closure of a maximal connected component of $V \setminus \aH$.  Since $(\aW,\aS)$ is irreducible, each alcove is a simplex.  We write $\fa$ for the alcove bounded by the hyperplanes $\{ H_{\alpha_1},\dots,H_{\alpha_n}, H_{\widetilde \alpha,1} \}$, and call $\fa$ the \emph{base alcove} or the \emph{fundamental alcove}.  The set of alcoves in $V$ is in bijection with the elements of $\aW$, and for $x \in \aW$ we write $\x$ or $x\fa$ for the alcove $x \cdot \fa$.  A \emph{panel} is a codimension one face of an alcove, and the \emph{supporting hyperplane} of a panel $p$ is the unique element of $\aH$ which contains $p$.

Let $\aH_0 = \{  H_{\alpha} \mid \alpha \in \Phi \}$, that is, $\aH_0$ is the set of all hyperplanes in $\aH$ which pass through the origin.  A \emph{Weyl chamber} is the closure of a maximal connected component of $V \setminus \aH_0$.  Each Weyl chamber is a simplicial cone, and we will also refer to Weyl chambers as sectors.  We define the \emph{dominant} Weyl chamber $\Cf$ to be the unique Weyl chamber which contains the base  alcove $\fa$. Thus $\Cf$ is the set of points $v \in V$ such that $\langle \alpha, v \rangle \geq 0$ for every $\alpha \in \Phi^+$. We denote by $\Cfm$ the \emph{antidominant} Weyl chamber, which is the unique chamber opposite $\Cf$ in $V$, consisting of all points $v \in V$ such that $\langle \alpha, v \rangle \leq 0$ for all $\alpha \in \Phi^+$. The set of Weyl chambers is in bijection with the elements of $\sW$, and for $w \in \sW$ we often write $\Cw_w$ for the Weyl chamber $w \cdot \Cf$.   Thus $\Cf = \Cw_\id$ and $\Cfm = \Cw_{w_0}$, where $\id$ is the identity element of $\sW$ and $w_0$ is its longest element.

For any root $\alpha\in\Phi$, any $k \in \Z$, and any $w\in\sW$, we denote by $\alpha^{k,w}$ the closed half-space of $V$ bounded by the hyperplane $H_{\alpha, k}$ that contains a subsector of the Weyl chamber $\Cw_w$.  In particular, for any $k \in \Z$ the half-space $\alpha^{k,\id}$ contains a subsector of the dominant Weyl chamber $\Cf$, and $\alpha^{k,w_0}$ contains a subsector of the antidominant Weyl chamber $\Cfm$.  The group $\aW$ acts naturally on the set $\{ \alpha^{k,w} \mid k \in \Z, w \in \sW \}$ via $(t^\mu u) \cdot \alpha^{k,w} = (u\alpha)^{k+\langle u\alpha, \mu\rangle, uw}$, where $\mu \in R^\vee$ and $u \in \sW$. 

Now let $X$ be an affine building of type $(W,S)$.  We regard $X$ as a simplicial complex, and for any apartment $\App$ of $X$ we may fix an identification of $\App$ with $V$ such that we may talk about hyperplanes, alcoves, Weyl chambers and so on in $\App$.  We will usually refer to the closed half-spaces of $\App$  determined by hyperplanes as \emph{half-apartments}.

\subsection{Galleries}\label{sec:galleries}

This section recalls definitions concerning combinatorial galleries that will be used in the sequel.  In this section, $X$ is any affine building of type $(\aW,\aS)$ irreducible, and $\App$ is any apartment of $X$.

We start with the following definition, which gives a special case of the combinatorial galleries of Gaussent and Littelmann (see Definition 8 of \cite{GaussentLittelmann}).

\begin{definition}\label{def:CombGallery}
	A \emph{(finite) combinatorial gallery} is a sequence of alcoves $c_i$ and faces $p_i$ in the affine building $X$ 
	$$ \gamma=(p_0, c_0, p_1, c_1, p_2, \dots , p_l, c_l, p_{l+1} ), $$
	where the first and last faces $p_0 \subseteq c_0$ and $p_{l+1}\subseteq c_l$ are simplices, and for $1 \leq i \leq l$ the face~$p_i$ is  a panel of both alcoves $c_{i-1}$ and $c_{i}$.    
\end{definition}

We will also need the (obvious) extension of this definition to infinite galleries.

\begin{definition}\label{def:InfiniteCombGallery}
	An \emph{infinite combinatorial gallery} is an (infinite) sequence of alcoves $c_i$ and faces $p_i$ in the affine building $X$ 
	$$ \gamma=(p_0, c_0, p_1 , c_1 , p_2 , \dots , p_i , c_i , p_{i+1} , \dots ), $$
	where the first face $p_0 \subseteq c_0$ is a simplex, and for all $i \geq 1$ the face $p_i$ is a panel of both alcoves $c_{i-1}$ and $c_{i}$.    
\end{definition}

We remark that if $c_{i-1} \neq c_{i}$ in either of these definitions, then there is no choice for the panel $p_i$.  We will often omit the word ``combinatorial" when referring to galleries.  All of our galleries will contain at least one alcove, and we will mostly be restricting attention to galleries that lie in a single apartment.

A gallery 
$ \gamma$ as in Definition~\ref{def:CombGallery} 
has \emph{length}~$l+1$, that is, its length is the number of alcoves counted with multiplicity.  We say that such a $\gamma$ is \emph{minimal} if it has minimal length among all galleries with first face $p_0$ and last face $p_{l+1}$.  An infinite gallery $\gamma$ as in Definition~\ref{def:InfiniteCombGallery} is \emph{minimal} if each of its initial finite subgalleries $(p_0, c_0, p_1 , c_1 , p_2 , \dots , p_i)$  is minimal.

Recall that each simplex in $X$ has a type which is a subset of $\aS$.  

\begin{definition}\label{def:type} Given a gallery $\gamma$ as in Definition~\ref{def:CombGallery}, the \emph{type} of $\gamma$ is given by the $(l+2)$-tuple 
\[
\type(\gamma) := (\type(p_0), s_{j_1}, s_{j_2}, \dots, s_{j_l}, \type(p_{l+1}))
\]
where for $1 \leq i \leq l$ the panel $p_i$ has type $s_{j_i} \in \aS$.  If $p_0 = c_1$ and $p_{l+1} = c_l$ then we will usually write instead $\type(\gamma) = (s_{j_1}, s_{j_2}, \dots, s_{j_l})$, so that the type of $\gamma$ is a word in $\aS$.
\end{definition}

It will not always be necessary to record the details of a gallery's type.  In particular, for $x \in \aW$ and a fixed minimal gallery from $\fa$ to $\x$, we will sometimes describe any gallery of this same type with first face $\fa$ as having \emph{type~$\vec{x}$}.  Similarly, for $\lambda \in R^\vee$ and a fixed minimal gallery from the origin to $\lambda$, we will sometimes describe any gallery of this same type with first face the origin as having \emph{type $\vec{\lambda}$}.

We finally recall several general definitions concerning orientations and positively folded galleries.  See the introduction and~\cite{GraeberSchwer} for some history concerning these notions.

\begin{definition}[Definition~3.1 of~\cite{GraeberSchwer}]\label{def:orientation} An \emph{orientation} of the apartment $\App$ is a map $\phi$ that assigns to every pair $(c,p)$ consisting of an alcove $c$ of $\App$ and a panel $p$ of $c$ one of the symbols $+$ or $-$. We say that an alcove $c$ is on the \emph{positive} side of $p$ if $\phi(c,p)=+$. 
\end{definition}

\begin{definition}\label{def:fold}
	Let $\gamma$ be a combinatorial gallery as in Definition~\ref{def:CombGallery}.  For $1 \leq i \leq l$, the gallery $\gamma$ is said to be \emph{folded at $p_i$}, or to \emph{have a fold at $p_i$}, if  $c_{i-1} = c_{i}$. 
\end{definition}	

\noindent Combinatorial galleries include not just the alcoves $c_i$ but also the panels $p_i$ precisely in order to record the panels in which folds (may) occur.  If $\gamma$ is folded at $p_i$ we may also say that \emph{$\gamma$ is folded at $H$}, where $H$ is the hyperplane supporting $p_i$.  
	
\begin{definition}[Definition~4.5 of~\cite{GraeberSchwer}]\label{def:positivelyFolded}
	Fix an orientation $\phi$ of $\App$ and suppose that a combinatorial gallery $\gamma$ lies in $\App$. We say that $\gamma$ is \emph{positively folded with respect to  $\phi$} if for all $i$ such that $\gamma$ is folded at $p_i$, one has $\phi(c_i, p_i)=+$.
\end{definition}


\section{Chimneys, retractions, orientations, and shadows}\label{sec:chimneysRetractionsGalleries}

Throughout this section, we use the notation of Section~\ref{sec:preliminaries}.  Thus $X$ is an arbitrary affine building of type $(\aW,\aS)$, with $(\aW,\aS)$ irreducible, and $\App$ is any apartment of~$X$, with a fixed identification of $\App$ with the vector space $V$.  We define chimneys for such $X$ in Section~\ref{sec:chimneys}.  In Section~\ref{sec:retractions} we formulate retractions from chimneys, and in Section~\ref{sec:orientations} we define the orientation induced by a chimney.  Section~\ref{sec:shadows} then defines shadows in our setting, and relates  shadows to retractions.

\subsection{Chimneys}\label{sec:chimneys}

In this section we give our definition of chimneys, in Definition~\ref{def:J-chimney} and its generalization Definition~\ref{def:Jy-chimney}.  We discuss the relationship between our definitions and related notions in the literature, and provide examples.

\begin{definition}\label{def:J-chimney}
Let $J$ be a subset of $\{ 1, \dots, n \}$, with corresponding root system $\Phi_J$.  The \emph{$J$-chimney (in $\App$)} is the following collection of half-apartments of $\App$:
\[
\xi_J\define 
 \left\{\alpha^{k,\id} \mid \alpha \in \Phi_J^+, k\leq 0 \right\} 
\cup
\left\{\alpha^{k, w_0} \mid \alpha \in \Phi_J^+, k\geq 1\right\}
\cup
\left\{\beta^{k, w_0} \mid \beta\in \Phi^+ \setminus \Phi_J^+, k\in \Z\right\}.
\] 
Moreover, for any collection of integers $\{n_\beta\in\Z \mid \beta\in \Phi^+ \setminus \Phi_J^+\}$, the corresponding \emph{$J$-sector} is the subcomplex of $\App$ given by 
\[
S_J\left(\{n_\beta\}_{\beta\in \Phi^+ \setminus \Phi_J^+}\right)
\define 
\left( \bigcap_{\alpha\in\Phi_J^+} \alpha^{0, \id} \cap\alpha^{1, w_0} \right) \cap
\left( \bigcap_{\beta\in \Phi^+ \setminus \Phi_J^+} \beta^{n_\beta, w_0}\right).
\]  
Write $S_{J}(0)$ for the $J$-sector defined by the all-zeros sequence, i.e. $n_\beta=0, \forall \beta\in \Phi^+ \setminus \Phi_J^+$.   
\end{definition}

Before we continue with some examples of chimneys and $J$-sectors, we make a couple of remarks on the origin and broader mathematical context of chimneys. 
 
\begin{remark}[Chimneys and bordifications]\label{rem:chimneyCL}
It is important to note that in our definition, a chimney is a collection of subsets of the apartment~$\App$, so a chimney is not itself a subset of $\App$.  
The intersection of all half-apartments in a chimney is in fact most often empty, and the union of all such half-apartments is most often $\App$.  One way to think of the $J$-chimney $\xi_J$ is to see it as a direction at infinity of the affine building~$X$ which is determined by the given filter of half-apartments. 
More formally, a $J$-chimney for nonempty $J$ is an alcove in the affine building associated to the Levi factor associated with $J$. In the rank $1$ case, that is, where $|J|=1$, this building is a panel tree. 

The $J$-sectors should be thought of as representatives within~$\App$ of the $J$-chimney, much like Weyl chambers (or sectors) in $\App$ are representatives of chambers in the spherical building at infinity. In other words, for a fixed $J$ the collection of all $J$-sectors can be thought of as the collection of subcomplexes of $\App$ which ``point towards" the $J$-chimney $\xi_J$, which lies at infinity.  

These ideas can be made precise as follows. One can show that $\xi_J$ describes a point in the combinatorial bordification $\mathscr{C}_{sph}(X)$ introduced by Caprace and L\'ecureux in Section~2.1 of~\cite{CapraceLecureux}, and that $J$-sectors are a special case of the generalized sectors introduced in Section~2.3 of~\cite{CapraceLecureux}.  We remark that Caprace and L\'ecureux are considering an arbitrary building $X$, and are generalizing results obtained for certain Bruhat--Tits buildings by Guivarc'h and R\'emy in~\cite{GuivarchRemy}.  Another related work is the (unpublished) PhD thesis of Charignon~\cite{Charignon}, who defines a compactification of an arbitrary locally finite affine building by a collection of affine buildings at infinity.  

In order to keep this work self-contained and focused on the case at hand, we decided to introduce our own definitions of chimneys and the associated sectors.  Although the work of Caprace and L\'ecureux~\cite{CapraceLecureux} takes place in the setting of an arbitrary building $X$, we restrict our attention to affine cases.  Our algebraic applications in Section~\ref{sec:galleriesCosets} are for affine buildings, and for $(W,S)$ affine, the associated spherical root system provides a convenient way to describe chimneys and $J$-sectors. In later proofs we will mostly work with $J$-sectors rather than the chimney itself, although we will need our definition of a chimney to formulate the definition of a chimney-induced orientation (see Definition~\ref{def:chimneyOrientation} below).  
\end{remark}

\begin{remark}[Historical origin of chimneys]\label{rem:chimneyRousseau}
The notion of a chimney in an affine building is not new.  Chimneys (in French, \emph{chemin\'ees}) were first introduced by Rousseau (see for example~\cite{Rousseau77} or~\cite{Rousseau01}).  In these works, a chimney is defined as the closure in the building of the union of a geodesic segment and a geodesic half-line which have the same origin and lie in the same apartment.   A formulation of chimneys equivalent to Rousseau's, and which relates chimneys to compactifications of affine buildings, appears in Section~4.2 of~\cite{Charignon}.  Our $J$-sectors are in fact a special case of Rousseau's chimneys.
\end{remark}

We now discuss some examples of $J$-chimneys and $J$-sectors. 

\begin{example}\label{eg:emptyset}
If $J = \emptyset$ then we have $\Phi_J =\emptyset$ and so the chimney $\xi_{\emptyset}$ is the collection of all half-apartments containing subsectors of the antidominant Weyl chamber $\Cw_{w_0}$.   The corresponding $J$-sectors are all of the translates of  $\Cw_{w_0}$, that is, the collection of all sectors in $\App$ whose intersection with $\Cw_{w_0}$ contains a subsector of $\Cw_{w_0}$.  In particular, the sector $S_{\emptyset}(0)$ is equal to the antidominant Weyl chamber $\Cw_{w_0}$.  Thus the chimney $\xi_{\emptyset}$ can be identified with the chamber at infinity represented by the antidominant Weyl chamber.
\end{example}

\begin{example}\label{eg:all}
If $J = \{ 1,\dots,n\}$ then $\Phi_J=\Phi$.  The $J$-chimney is the collection of all half-apartments containing the base alcove $\fa$, and the only $J$-sector is the base alcove $\fa$ itself. 
\end{example}	

\begin{example}\label{eg:intro}  Consider Figures~\ref{fig:chimneyIntro} and~\ref{fig:AlcoveShadowA2} of the introduction.  Letting $\alpha = \alpha_1$ be the black root in these figures, the other two (gray) roots can be labeled as $\alpha_2$ and $\alpha_1 + \alpha_2$, so that the three depicted roots are the positive roots in type $\tilde{A_2}$.  Then if $J = \{ 1 \}$, the shaded region is an example of a $J$-sector.  
\end{example}

Recall that the affine Weyl group $\aW$ acts naturally on the collection of half-apartments in~$\App$.  We now use this action to generalize $J$-chimneys and $J$-sectors.

\begin{definition}\label{def:Jy-chimney}
For any $J \subseteq \{ 1,\dots,n\}$ and any $y \in \aW$, the \emph{$(J,y)$-chimney} $\xi_{J,y}$ is obtained by acting on the $J$-chimney on the left by $y$, that is, \[ \xi_{J,y}\define y \cdot \xi_J.\]  Hence, letting $y = t^\mu u$ with $\mu \in R^\vee$ and $u \in \sW$, 
\begin{eqnarray*}
\xi_{J,y} & \define & 
 \left\{(u\alpha)^{k+\langle u\alpha, \mu\rangle,u} \mid \alpha\in\Phi_J^+, k\leq 0 \right\}
\cup\left\{(u\alpha)^{k+\langle u\alpha, \mu\rangle, uw_0} \mid \alpha\in\Phi_J^+, k\geq 1\right\} \\
&& \cup
\left\{(u\beta)^{k, uw_0} \mid \beta\in \Phi^+ \setminus \Phi_J^+, k\in \Z\right\}.
\end{eqnarray*}
Similarly we obtain a \emph{$(J,y)$-sector} 
$ S_{J,y}(\{n_\beta\})\define y \cdot S_{J}(\{n_\beta\})$
for all $y=t^\mu u\in\aW$ and any set $\{n_\beta\in\Z \mid \beta\in \Phi^+ \setminus \Phi_J^+\}$ by 
\[
S_{J,y}\left(\{n_\beta\}_{\beta\in \Phi^+ \setminus \Phi_J^+}\right)
\define 
\left( \bigcap_{\alpha\in\Phi_J^+} (u\alpha)^{\langle u\alpha, \mu\rangle, u} \cap(u\alpha)^{1+\langle u\alpha, \mu\rangle, uw_0} \right) \cap
\left( \bigcap_{\beta\in \Phi^+ \setminus \Phi_J^+} (u\beta)^{n_\beta+\langle u\beta, \mu\rangle, uw_0}\right).
\]  
Write $S_{J,y}(0)$ for the $(J,y)$-sector defined by the all-zeros sequence, i.e. $n_\beta=0, \forall \beta\in \Phi^+ \setminus \Phi_J^+$.   
\end{definition}

The next two examples generalize Examples~\ref{eg:emptyset} and~\ref{eg:all}, respectively.

\begin{example}\label{eg:w-emptyset}
If $J = \emptyset$ and $w \in \sW$, then any $(J,w)$-sector is a translate of the Weyl chamber $w\Cw_{w_0} = \Cw_{ww_0}$, and the sector $S_{\emptyset,w}(0) = w \cdot S_\emptyset(0)$ is equal to the Weyl chamber $\Cw_{ww_0}$.  Thus the $(\emptyset,w)$-chimney can be identified with the chamber at infinity represented by the Weyl chamber $\Cw_{ww_0}$, and the chambers at infinity (in the boundary of~$\App$) are in bijection with the $(\emptyset,w)$-chimneys.
\end{example}

\begin{example}\label{eg:y-all}
If $J = \{1,\dots,n\}$ and $y \in \aW$ then the only $(J,y)$-sector is the alcove $y\fa$.  Thus every alcove in $\App$ is an instance of a $(J,y)$-sector.
\end{example}	

\subsection{Retractions from chimneys}\label{sec:retractions}

In this section we give our formulation of chimney retractions (see Definition~\ref{def:chimney-retraction}).  In order to give this definition and show that chimney retractions are well-defined, we will need several preliminary results and definitions.  Throughout this section, $J \subseteq \{1,\dots,n\}$ and $y \in \aW$.
 
\begin{lemma}\label{lem:one-alcove}  Every $(J,y)$-sector contains at least one alcove, and if it contains more than one alcove, it contains infinitely many.
\end{lemma}
\begin{proof}
It is enough to prove the statement for $J$-sectors.

Suppose first that $J = \emptyset$.  Then by Example~\ref{eg:emptyset}, any $J$-sector contains a subsector of the antidominant Weyl chamber $\mathcal{C}_{w_0}$, hence contains infinitely many alcoves. 

Now suppose $J = \{ 1,\dots,n\}$.  Then by Example~\ref{eg:all}, the only $J$-sector is the base alcove~$\fa$, so the statement holds in this case. 

In all other cases, the sets $\Phi_J^+$ and $\Phi^+ \setminus \Phi^+_J$ are both nonempty, and any $J$-sector is contained in the region of $\App$ lying between $H_\alpha$ and $H_{\alpha,1}$ for all $\alpha \in \Phi_J^+$.  The conditions provided by the (finite) collection of roots in $\Phi^+ \setminus \Phi^+_J$ then determine a cone containing infinitely many alcoves of this region of $\App$.  Hence the assertion of the Lemma. 
\end{proof}

\begin{definition}
Let $y=t^\mu u\in \aW$, where $\mu \in R^\vee$ and $u \in \sW$.  Let $S_{J,y}(\{n_\beta\})$ be a $(J,y)$-sector.  A minimal infinite gallery 
$$ \gamma=(p_0, c_0, p_1 , c_1 , p_2 , \dots , p_i , c_i , p_{i+1} , \dots )$$
is \emph{regular} with respect to $S_{J,y}(\{n_\beta\})$ if for all $\beta\in \Phi^+ \setminus\Phi_J^+$ and every $k \in \Z$ such that the half-apartment $(u\beta)^{k, uw_0}$ is a subset of the half-apartment $y \cdot \beta^{n_\beta,w_0} = (u\beta)^{n_\beta+\langle u\beta, \mu\rangle, uw_0}$, there exists an index $i_k$ such that the alcove $c_{i_k}$ of $\gamma$ is contained in $(u\beta)^{k, uw_0}$. 
\end{definition}

The previous definition essentially means that a regular infinite gallery eventually crosses every hyperplane that ``cuts across" the $(J,y)$-sector, and can only be trapped in a ``strip" of the $(J,y)$-sector lying between two parallel hyperplanes if they are of type $u\alpha$ with $\alpha \in \Phi_J^+$. 

For the next lemma we use the metric approach to buildings, details of which can be found in Section 12 of \cite{AB}.  Thus we view the affine building $X$ as a $\CAT(0)$ space equipped with the usual $\CAT(0)$ metric, which restricts to a Euclidean metric on each apartment.

\begin{lemma}\label{lem:gallery-support}
For every pair of alcoves $c$ and $d$ in $X$, and any points $x\in c$ and $y\in d$, there exists a minimal gallery containing the (unique) geodesic segment from $x$ to $y$. Moreover, every geodesic ray starting at $x$ is contained in a minimal infinite gallery. 
\end{lemma}
\begin{proof}
The statement about pairs of alcoves is taken from \cite[Lemma 3.1]{Marquis}. To see the statement about infinite geodesic rays argue as follows.  Each geodesic ray is contained in an apartment. By the first part of the current lemma each finite length initial piece is contained in a minimal gallery. Local finiteness of the Coxeter complex then implies that we may choose these initial pieces in a consistent way such that the resulting infinite gallery covering the geodesic ray is minimal. 
\end{proof}

\begin{corollary}\label{cor:regular-gallery}
Every $(J,y)$-sector that is not a single alcove contains a regular minimal infinite gallery. 
\end{corollary}
\begin{proof} 
A $(J,y)$-sector is by definition an intersection of half-apartments, hence $(J,y)$-sectors are convex, meaning that any two alcoves in a $(J,y)$-sector can be connected by a minimal gallery which lies in that sector.  Now by Lemma~\ref{lem:one-alcove}, a $(J,y)$-sector which is not a single alcove contains infinitely many alcoves.  By the proof of Lemma~\ref{lem:one-alcove} and the specific shape of the $(J,y)$-sector described there, a $(J,y)$-sector which is not a single alcove contains a regular geodesic ray starting at a point $x$ contained in the interior of an alcove $c$. This ray is, by Lemma~\ref{lem:gallery-support}, contained in a minimal infinite gallery $\gamma$ beginning at $c$. Moreover $\gamma$ is regular since the geodesic was regular.  By convexity of the intersection of half-apartments, $\gamma$ must be contained in the $(J,y)$-sector. 
\end{proof}

We will use the following result of~\cite{Charignon}, which holds for arbitrary affine buildings (in fact for arbitrary buildings).  Recall that any building has a complete apartment system.

\begin{lemma}[Corollary 4.3.2 in \cite{Charignon}]\label{lem:subgallery}
For every alcove $c$ in the affine building $X$ and any minimal infinite gallery $\gamma$ in $X$, there exists an apartment $\App'$ in the complete apartment system of $X$ that contains both $c$ and an infinite subgallery of $\gamma$. 
\end{lemma}

The next result will be crucial for our definition of a chimney retraction. 

\begin{prop}\label{prop:J-apartment}
Let $J \subseteq \{1,\dots,n\}$ and let $y \in \aW$.  Then for every alcove $c$ in $X$ there exists a collection of integers $\{n_\beta\in\Z \mid \beta\in \Phi^+ \setminus\Phi_J^+\}$  and an apartment $\App_{c,(J,y)}$ in the complete apartment system of $X$ such that $\App_{c,(J,y)}$ contains both $c$ and the $(J,y)$-sector $S_{J,y}(\{n_\beta\})$.  
\end{prop}
\begin{proof}
In the case that the unique $(J,y)$-sector is the alcove $y\fa$ this follows directly from the buildings axioms. 
In all other cases the set $\Phi^+ \setminus \Phi^+_J$ is nonempty, and so $S_{J,y}(0)$ contains by Corollary~\ref{cor:regular-gallery} a regular minimal infinite gallery $\gamma$.  Then Lemma~\ref{lem:subgallery} yields an apartment $\App'$ which contains both $c$ and an infinite subgallery of $\gamma$.  Put $\App_{c,(J,y)} = \App'$.  By the regularity of $\gamma$ and the fact that apartments intersect in convex sets, one obtains that the intersection of $\App_{c,(J,y)}$ with the apartment $\App$ does contain a $(J,y)$-sector $S_{J,y}(\{n_\beta\})$ for a collection of integers $n_\beta$ large enough.   	
\end{proof}

Proposition~\ref{prop:J-apartment} is the natural generalization of the axiomatic fact that any pair of alcoves is contained in a common apartment, and at the same time a natural generalization of the statement of part (1) of the Theorem in \cite[Section 8]{Brown} to chimneys. Corollary~\ref{cor:union} below is a direct consequence of the aforementioned proposition. 
 
\begin{corollary}\label{cor:union}
For any $J \subseteq \{1,\dots,n\}$, any $y \in \aW$, and any $(J,y)$-chimney, the building $X$ is (as a set) the union of all apartments containing a $(J,y)$-sector representing this chimney.  
\end{corollary}

We now combine Proposition~\ref{prop:J-apartment} with Lemma~\ref{lem:one-alcove} to obtain the following. 

\begin{corollary}\label{cor:intersection}
For any choice of an apartment  $\App_{c,(J,y)}$ as in the statement of Proposition~\ref{prop:J-apartment}, the intersection $\App \cap \App_{c,(J,y)}$ contains at least one alcove.  Hence there is a unique isomorphism $\App_{c,(J,y)} \to \App$ which fixes $\App_{c,(J,y)} \cap \App$ pointwise.
\end{corollary}

Proposition~\ref{prop:J-apartment} and Corollary~\ref{cor:intersection} are the reasons why chimney retractions as formulated in the next definition are well-defined.  We note that we are using an apartment $\App_{c,(J,y)}$ in the complete apartment system on $X$ to define this retraction.  
		
\begin{definition}\label{def:chimney-retraction}
Let $X$ be an affine building and $\App$ an apartment of $X$.  Let $J \subseteq \{1,\dots,n\}$ and let $y \in \aW$.   
For any alcove $c$ of $X$, let $r_{J,y}(c)$ be the image of $c$ under the unique isomorphism that maps an apartment $\App_{c,(J,y)}$ as in the statement of Proposition~\ref{prop:J-apartment} onto $\App$ while fixing $\App_{c,(J,y)} \cap \App$ pointwise.  The resulting induced simplicial map  $$r_{J,y}:X\to \App$$ is the \emph{$(J,y)$-chimney retraction} or the \emph{retraction from the $(J,y)$-chimney}.  
\end{definition}

A \emph{chimney retraction} or a \emph{retraction from a chimney} is a $(J,y)$-chimney retraction for some $J \subseteq \{ 1,\dots,n\}$ and some $y \in \aW$.  If $y = \id$ is the trivial element of $\aW$, we may write $r_J = r_{J,\id}$, and call $r_J$ the \emph{retraction from the $J$-chimney}.

\begin{remark}\label{rem:retractionClassical}  If $J = \emptyset$ and $y = w \in \sW$, so that the $(J,w)$-chimney can be identified with a chamber at infinity (see Example~\ref{eg:w-emptyset}), then the retraction from the $(J,w)$-chimney is the retraction onto $\App$ centered at this chamber at infinity.  If $J = \{1,\dots,n\}$, so that the unique $(J,y)$-sector is the alcove $y\fa$ (see Example~\ref{eg:y-all}), then the retraction from the $(J,y)$-chimney is the retraction onto $\App$ centered at the alcove $y\fa$.   Thus in these cases the chimney retraction specializes to the two well-known retractions of affine buildings onto an apartment.   
\end{remark}

\begin{remark}\label{rem:GHKR}
For arbitrary $J$, our formulation of the retraction $r_{J,y}$ is a generalization and formalization of the retraction $\rho_{P,w}$ which is described in Section~11.2 of~\cite{GHKRadlvs}, and also considered in Section~6 of~\cite{HKM}.  We explain this further at Remark~\ref{rem:GHKR-retraction}.  
\end{remark}

\begin{remark}\label{rem:retractionNonaffine}
A slightly different approach to chimney retractions, which was explained to us by Jean L\'ecureux, would be to use Proposition~2.30 of~\cite{CapraceLecureux}, or similar results from~\cite{Charignon}.  Since the definitions and results of~\cite{CapraceLecureux} hold for arbitrary buildings $X$, chimney retractions could then be defined for non-affine buildings as well (see also Remark~\ref{rem:chimneyCL}).
\end{remark}

The next corollary is an immediate consequence of our formulation of chimney retractions, and generalizes the corresponding result for retractions centered at chambers at infinity.  See the second paragraph after Exercise 11.65 of~\cite{AB}.

\begin{corollary}\label{cor:stable-alcove}  Let $X$ be an affine building and $\App$ an apartment of $X$.  Let $J \subseteq \{1,\dots,n\}$ and let $y \in \aW$.  Then for any alcove $c$ of $X$, there is an alcove $d$ in some $(J,y)$-sector of $\App$ such that the image of the alcove $c$ under the retraction $r_{J,y}: X \to \App$ is equal to the image of $c$ under the retraction onto $\App$ centered at $d$.
\end{corollary}

The following generalizes Corollary~\ref{cor:stable-alcove}, and extends to all chimney retractions the corresponding result for retractions centered at chambers at infinity (see Exercise 11.67 of~\cite{AB}).  As explained further at Remark~\ref{rem:GHKR-retraction}, this result also generalizes Lemma~6.4 of~\cite{HKM}.

\begin{corollary}\label{cor:sim-stable-alcove}  Let $X$ be an affine building and $\App$ an apartment of $X$.  Let $J \subseteq \{1,\dots,n\}$ and let $y \in \aW$.  Then for any bounded subset $B$ of $X$, there is an alcove $d$ in some $(J,y)$-sector of $\App$ such that the image of $B$ under the retraction $r_{J,y}: X \to \App$ is equal to the image of $B$ under the retraction onto $\App$ centered at $d$.
\end{corollary}

\begin{proof}  If $J= \{1,\dots,n\}$ then the unique $(J,y)$-sector is the alcove $\y = y\fa$, so taking $d = \y$ we are done.  Now assume that $J \subsetneq \{1,\dots,n\}$ and fix a $(J,y)$-sector in $\App$.  Then by Corollary~\ref{cor:regular-gallery}, this $(J,y)$-sector contains a regular infinite gallery $ \gamma=(p_0, c_0, p_1, \dots , p_i , c_i , p_{i+1} , \dots )$.  By similar arguments to those used to prove Proposition~\ref{prop:J-apartment}, since $B$ is bounded, for all $i$ large enough (depending on $B$) the restriction to $B$ of $r_{J,y}: X \to \App$ is equal to the restriction to $B$ of the retraction onto $\App$ centered at $c_i$.  This completes the proof.
\end{proof}

\subsection{Chimney-induced orientations and positively folded galleries}\label{sec:orientations}

In this section we define the natural orientation induced by a $(J,y)$-chimney, and hence say what it means for a combinatorial gallery to be positively folded with respect to a given chimney.   The general definitions of an orientation and a positively folded gallery are stated in Section~\ref{sec:galleries}.

The key definition for this section is as follows. Recall from Definition~\ref{def:Jy-chimney} that a $(J,y)$-chimney is a collection of half-apartments in $\App$.

\begin{definition}\label{def:chimneyOrientation} 
Let $J \subseteq \{1,\dots,n\}$ and let $y \in \aW$.  Let $c$ be an alcove of $\App$, let $p$ be a panel of $c$, and let $H$ be the hyperplane supporting $p$.  The \emph{orientation induced by the $(J,y)$-chimney}, denoted $\phi_{J,y}$, is defined by $\phi_{J,y}(c,p) = +$ if and only if the half-apartment determined by $H$ which contains $c$ is \emph{not} an element of the $(J,y)$-chimney.  
\end{definition}

\noindent Thus an alcove $c$ is on the positive side of $p$ with respect to the orientation $\phi_{J,y}$ if, roughly speaking, $c$ is on the side of $p$ which faces away from the $(J,y)$-chimney.  We say that a gallery in $\App$ is \emph{positively folded with respect to the $(J,y)$-chimney} if it is positively folded with respect to the orientation $\phi_{J,y}$ (see Definition~\ref{def:positivelyFolded}).  Thus a gallery is positively folded with respect to a given chimney if, roughly speaking, it is always folded away from this chimney.

The next result is immediate from Definition~\ref{def:chimneyOrientation}.  

\begin{lemma} \label{lem:side} Let $\phi = \phi_{J,y}$ be the orientation induced by the $(J,y)$-chimney and let $H$ be a hyperplane of $A$.  If $p$ and $p'$ are panels of $H$, and $c$ and $c'$ are alcoves of $A$ containing $p$ and~$p'$, respectively, then $\phi(c,p) = \phi(c',p')$ if and only if $c$ and $c'$ are on the same side of $H$ in $A$.
\end{lemma}

\noindent We now use Lemma~\ref{lem:side} to make the following definition for chimney orientations (for an arbitrary orientation, it does not make sense to talk about \emph{the} positive or negative side of a hyperplane). 

\begin{definition}\label{def:positive-negative}  Let $\phi = \phi_{J,y}$ be the orientation induced by the $(J,y)$-chimney.  Let $H$ be a hyperplane of $A$, let $p$ be a panel in $H$ and let $c$ be an alcove containing $p$.  Then $c$ is on \emph{the positive side of $H$ (with respect to $\phi$)} if $\phi(c,p) = +$, and $c$ is on \emph{the negative side of $H$ (with respect to $\phi$)} if $\phi(c,p) = -$.\end{definition}

\begin{remark}\label{rem:orientation}
The orientation induced by a chimney generalizes several notions of orientation that already appeared in the literature.

If $J = \emptyset$ then using Example~\ref{eg:emptyset}, one can check that the orientation induced by the $J$-chimney is the same as the periodic orientation on hyperplanes considered in~\cite{PRS}.  More generally, for $w \in \sW$, the orientation induced by the  $(J,w)$-chimney is the same as the periodic orientation on hyperplanes induced by the labeling at infinity $\phi^\partial_w$ as in Definition~3.7 of~\cite{MST1} (see also Definitions~3.10 and~3.13 of~\cite{GraeberSchwer}); such an orientation is often called an orientation at infinity in~\cite{MST1}.  In particular, the opposite standard orientation at infinity in~\cite{MST1} is that induced by  the $\emptyset$-chimney, and this orientation is what is used for the main constructions in~\cite{MST1}.

If $J = \{1,\dots,n\}$ and $y \in \aW$ then the orientation induced by the $(J,y)$-chimney is opposite to the alcove orientation towards $y\fa$ from Definition 3.6 of~\cite{GraeberSchwer}.  More precisely, recall from Example~\ref{eg:all} that $\y = y\fa$ is the unique $(J,y)$-sector, and write $\phi_\y$ for the alcove orientation towards $y\fa$.  Then $\phi_{J,y}(c,p) = +$ if and only if $\phi_\y(c,p) = -$.
\end{remark}

We now discuss braid-invariance of orientations.  This will be used in Section~\ref{sec:shadows} below.  

\begin{definition}[{\cite[Definition 5.2]{GraeberSchwer}}]\label{def:braid invariance} 
An orientation $\phi$ of $\App$ is \emph{braid-invariant} if for any braid-equivalent reduced words $(s_{j_1},\dots,s_{j_l})$ and $(s'_{j_1},\dots,s'_{j_l})$ in $\aS$ and any $x \in \aW$, there is a combinatorial gallery of type $(\type(p_0),s_{j_1},\dots,s_{j_l},\type(p_{l+1}))$ with final alcove $\x$ which is positively folded with respect to $\phi$ if and only if there is a combinatorial gallery of type $(\type(p_0),s'_{j_1},\dots,s'_{j_l},\type(p_{l+1}))$ with final alcove $\x$ which is positively folded with respect to $\phi$. 
\end{definition}

It is proved in Proposition~4.33 of~\cite{MST1} (see also Proposition~5.3 of~\cite{GraeberSchwer}) that for $J = \emptyset$ and any $w \in \sW$, the orientation induced by the $(J,w)$-chimney is braid-invariant.  (Example~\ref{eg:emptyset} and Remark~\ref{rem:orientation} explain the relationship between this language and the terminology used in~\cite{MST1} and~\cite{GraeberSchwer}.)  We now establish braid-invariance for arbitrary chimneys.  We will describe in Remarks~\ref{rem:braidRecursion} and~\ref{rem:braidDoubleCoset} how our main results give two other approaches to the proof of the following. 

\begin{prop}\label{prop:braidinvariant}
	Chimney-induced orientations are braid-invariant. 
\end{prop}	

\begin{proof}  Let $\phi = \phi_{J,y}$ be a chimney-induced orientation.  As in the proofs of \cite[Proposition~4.33]{MST1} and~\cite[Proposition~5.3]{GraeberSchwer}, it suffices to consider reduced words $(s,t,s,\dots)$ and $(t,s,t,\dots)$ of length $m$, where $s, t \in S$ are such that $st$ has order $m = m_{st} < \infty$, and galleries of these types with first alcove $\fa$.  All panels of such galleries are panels of types $s$ or $t$ which contain the origin, and all alcoves of such galleries are of the form $v\fa$ where $v \in \langle s, t \rangle$.  It is then straightforward to verify that the orientation induced by the $(J,y)$-chimney is ``locally" the same as the orientation induced by the $(\emptyset,w)$-chimney, for some $w \in \sW$ (depending on $s$, $t$, $J$, and $y$).  To be precise, one can check that there is a $w \in \sW$ such that for every pair $(c,p)$ where $c = v\fa$ with $v \in \langle s, t \rangle$, and $p$ is a panel of $c$ of type $s$ or $t$, we have $\phi(c,p) = \phi_{\emptyset,w}(c,p)$.  The result then follows from~\cite[Proposition~4.33]{MST1} or~\cite[Proposition~5.3]{GraeberSchwer}.
\end{proof}

In the rest of the paper we will only consider orientations $\phi=\phi_{J,y}$ which are induced by some $(J,y)$-chimney.

\subsection{Shadows}\label{sec:shadows}

In this section we recall and slightly generalize the notion of a shadow which was first introduced in \cite{GraeberSchwer}. 
A special case of Definition~\ref{def:shadow} below, namely that where $\sigma = \tau = \fa$, is the notion of a regular shadow as introduced in \cite[Definition 6.3]{GraeberSchwer}. Proposition~\ref{prop:braidinvariant} implies that shadows are well-defined and independent of the choice of minimal gallery $\gamma$ in the following definition. 

\begin{definition}\label{def:shadow}  Let $J \subseteq \{1,\dots,n\}$ and $y \in \aW$. Let $x \in \aW$ and let $\sigma$ and $\tau$ be faces of the fundamental alcove $\fa$ which contain the origin.  Fix a minimal gallery $\gamma$ with first face $\sigma$ and last face $x\tau$. 
The \emph{shadow of $x\tau$ starting at $\sigma$ with respect to the $(J,y)$-orientation}, denoted $\Sh_{J,y}(x\tau, \sigma)$, is the set of final simplices of all galleries of $\type(\gamma)$ which have first face~$\sigma$ and are positively folded with respect to the $(J,y)$-chimney. 
\end{definition}

We sometimes say that a shadow is taken \emph{with respect to the $(J,y)$-chimney} instead of with respect to the (chimney-induced) $(J,y)$-orientation.  

\begin{notation}\label{not:shadow}  We simplify terminology and notation for the most commonly used cases.

If $J=\{1,\ldots,n\}$, then the $(J,y)$-chimney is just the alcove $\y=y\fa$. In this case, we will denote the shadow by $\Sh_{\y}(x\tau,\sigma).$

	If both $\sigma$ and $\tau$ are equal to the fundamental alcove $\fa$, so that $x\tau = x\fa = \x$, the \emph{shadow of $\x$ with respect to the $(J,y)$-chimney} is defined to be the shadow of $\x$ starting at $\fa$ with respect to the $(J,y)$-orientation.  That is, $\Sh_{J,y}(\x)$ is the set of final alcoves of galleries of type $\vec{x}$ which have first face and first alcove $\fa$ and which are positively folded with respect to the $(J,y)$-chimney.   We denote this shadow by $\Sh_{J,y}(\x)$.  If in addition $J = \{1,\dots,n\}$, we denote this shadow by $\Sh_\y(\x)$.
	
	If $\sigma$ and $\tau$ are both the origin, then $\lambda:=x\tau$ is an element of the coroot lattice.  In this case, the \emph{shadow of $\lambda$ with respect to the $(J,y)$-chimney} is the set of final vertices of galleries of type $\vec{\lambda}$ which have first face the origin and which are positively folded with respect to the $(J,y)$-chimney.  We denote this shadow by $\Sh_{J,y}(\lambda)$.  If in addition $J = \{1,\dots,n\}$, we denote this shadow by $\Sh_\y(\lambda)$.
\end{notation}

Examples of shadows are given in Figures~\ref{fig:chimneyIntro} and~\ref{fig:AlcoveShadowA2} of the introduction, and in Section~\ref{sec:examples}.

The next result gives an interpretation of shadows in terms of retractions.  We write $r:X \to \App$ for the retraction onto $\App$ centered at the base alcove~$\fa$.  Recall from Definition~\ref{def:chimney-retraction} that $r_{J,y}: X\to \App$ denotes the retraction (onto $\App$) from the $(J,y)$-chimney.  

\begin{prop}\label{prop:shadowsRetractions}
Let $J \subseteq \{1,\dots,n\}$ and let $y \in \aW$.      
 
Let $x \in \aW$ and let $\sigma$ and $\tau$ be faces of the fundamental alcove $\fa$ which contain the origin.  Let $W_\sigma$ be the subgroup of $\sW$ which fixes $\sigma$.  Then the shadow of $x\tau$ with respect to the $(J,y)$-orientation starting at $\sigma$ is the set of simplices given by 
	\[
	\Sh_{J,y}(x\tau,\sigma) = r_{J,y} \circ r^{-1}(W_\sigma \cdot x\tau).
	\]
In particular, the shadow of $\x = x\fa$ with respect to the $(J,y)$-chimney is the set of alcoves given by 	
	\[
	\Sh_{J,y}(\x) =  r_{J,y} \circ r^{-1}(\x).
	\]

For $\lambda \in R^\vee$, the shadow of $\lambda$ with respect to the $(J,y)$-chimney is the set of elements of $R^\vee$ given by 
	\[
	\Sh_{J,y}(\lambda) = r_{J,y} \circ r^{-1}(\sW \cdot \lambda).
	\]
\end{prop}

\begin{proof}
Fix a minimal gallery $\gamma$ with first face  $\sigma$ and last face $x\tau$. By definition the shadow $\Sh_{J,y}(x\tau,\sigma)$ is then the collection of end-simplices of all galleries of $\type(\gamma)$ and first face $\sigma$. 	

As described in Section~4.3 of~\cite{GraeberSchwer}, the affine Weyl group $\aW$ acts naturally on the set of combinatorial galleries in $\App$.  This action is type-preserving and also preserves minimality of galleries.  Thus the stabilizer $\aW_\sigma$ of the first face $\sigma$ inside $\aW$ permutes the set of all minimal galleries of $\type(\gamma)$ in $\App$ which have first face $\sigma$, and $\aW_\sigma \cdot x\tau$ is the collection of all end-simplices of minimal galleries in $\App$ which start in $\sigma$ and are of $\type(\gamma)$.  Hence by definition of the retraction $r$, the set $r^{-1}(\aW_\sigma\cdot x\tau)$ is the collection of all end-simplices of minimal galleries in $X$ which start in $\sigma$ and are of $\type(\gamma)$.

We prove first that the image of every such minimal gallery under $r_{J,y}$ is positively folded with respect to the $(J,y)$-chimney, and has $\type(\gamma)$.  The proof is exactly the same as the proof of Lemma~2.9 in~\cite{Hitzel}, where one replaces the sector $\mathcal{C}_f^{op}$ by a $(J,y)$-sector. The main argument then makes use of Proposition~\ref{prop:J-apartment} and Corollary~\ref{cor:stable-alcove}.

It then remains to show that every gallery in $\App$ of $\type(\gamma)$ which starts at $\sigma$ and is positively folded with respect to the $(J,y)$-chimney is an image of a minimal one under the retraction $r_{J,y}$. This is the analog of Proposition~3.3 in~\cite{Hitzel} for  chimney retractions. The proof of this proposition also carries over to the chimney case if we prove the analog of \cite[Lemma 2.4]{Hitzel}. This is done in Lemma~\ref{lem:compatability} below.
\end{proof}

The next lemma completes the last step in the proof of Proposition~\ref{prop:shadowsRetractions}.  Note that in the following statement, the $(J,y)$-sectors are each contained in our fixed apartment $\App$, but that once any apartment $A_i$ contains some $(J,y)$-sector, we can also define the retraction onto $A_i$ from the $(J,y)$-chimney.

\begin{lemma}\label{lem:compatability}  Let $J \subseteq \{1,\dots,n\}$ and $y \in \aW$.
Let $A_1,A_2,\dots,A_k$ be a collection of apartments in $X$ such that each $A_i$ contains some $(J,y)$-sector. Denote by $r_i$ the retraction onto~$A_i$ from the $(J,y)$-chimney. Then 
\[r_1\circ r_2 \circ \cdots \circ r_k = r_1.\]  
\end{lemma}
\begin{proof}
For all $i$ the retraction $r_i$ maps, by definition, any apartment containing a $(J,y)$-sector isomorphically onto $A_i$. Corollary~\ref{cor:union} implies that $X$ is the union of such apartments. 
Define $\rho\define r_1\circ\cdots\circ r_k$. Then on the one hand each apartment that contains a $(J,y)$-sector is isomorphically mapped onto $A_1$ by $\rho$. On the other hand $\rho$ is (Weyl-)distance nonincreasing. These two properties, which also hold for $r_1$, uniquely characterize this retraction. Therefore $r_1=\rho$.  
\end{proof}

\section{A recursive description of shadows}\label{sec:recursion}

The main result in this section is Theorem~\ref{thm:AlcoveRecursion}, which gives a recursive description of shadows.  A special case of this theorem was stated as Theorem~\ref{thm:recursionIntro} in the introduction.  The examples in Section~\ref{sec:examples} illustrate how to apply our recursion.
  
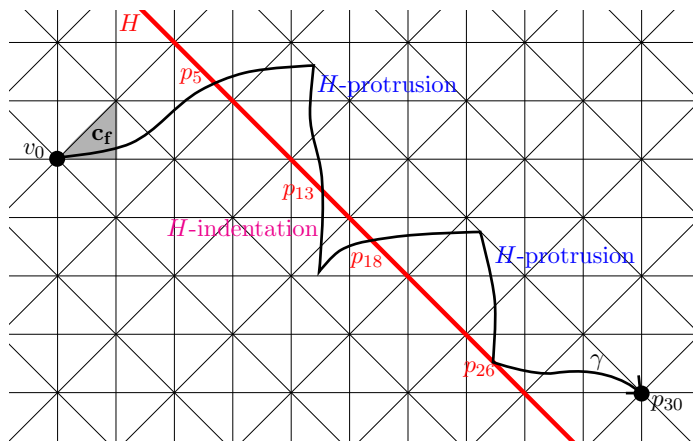
\begin{figure}[ht]
	\begin{center}
		\resizebox{0.6\textwidth}{!}
		{
			\begin{overpic}{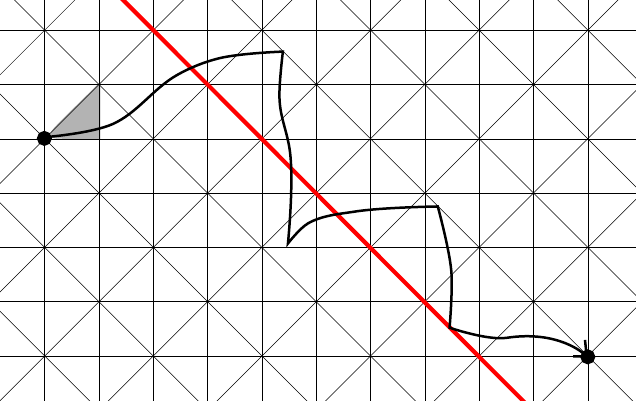}
				\put(16,60){\color{red}{$H$}}
				\put(85,11){$\gamma$}
				\put(2,42){\color{black}{$v_0$}}
				\put(12,44){\color{black}{$\fa$}}
				\put(25,53){\color{red}{$p_5$}}
				\put(40,36){\color{red}{$p_{13}$}}
				\put(50,26){\color{red}{$p_{18}$}}
				\put(66.5,10){\color{red}{$p_{26}$}}
				\put(94,5){\color{black}{$p_{30}$}}
				\put(45,51){\color{blue}{$H$-protrusion}}
				\put(71,26){\color{blue}{$H$-protrusion}}
				\put(23,30){\color{magenta}{$H$-indentation}}
			\end{overpic}
		}
		\caption{A gallery $\gamma$ with $H$-protrusions and an $H$-indentation. See Example~\ref{Ex:OutcropsIngrowths} for further details. }
		\label{fig:OutcropsIngrowths}
	\end{center}
\end{figure}

We continue all notation from Section~\ref{sec:chimneysRetractionsGalleries}.  Let $H$ be a hyperplane in $A$.  For $y \in \aW$, we denote by $H^y$ the (closed) half-apartment of $A$ which is bounded by $H$ and contains the alcove $\y = y\fa $.  In particular, $H^\id$ is the closed half-apartment of $A$ which is bounded by $H$ and contains the base alcove $\fa = \id\fa$.  For example, in Figure~\ref{fig:OutcropsIngrowths}, $H^\id$ is the half-apartment to the left of $H$.   We write $r_H$ for the reflection of $A$ in $H$, and note that $r_H$ is type-preserving.

We now define two operators on galleries, $e^\L_H$ and $f^{\L}_H$.  We will use these together with Lemma \ref{lem:GammaFold} to prove Theorem~\ref{thm:AlcoveRecursion}.  Some parts of the following definitions are illustrated by Figure~\ref{fig:OutcropsIngrowths}.

\begin{definition} \label{def:OutcropsIngrowths}
Fix a hyperplane $H$ of $A$, and let 
$\gamma=(p_0, c_0 , p_1, c_1 , p_2 , \ldots , p_l , c_l ,  p_{l+1})$
be a combinatorial gallery in $A$.
\begin{enumerate}[label={(\roman*)}]
	\item For $0 \leq j < k \leq l+1$, call ${\L}=[j,k] \subset \N$ an \emph{$H$-protrusion (of $\gamma$)} if $p_j,p_k\subset H$ and $c_i\not\subset H^\id$ for all $j\leq i <k$ (note that $p_i \subset H$ may occur for some $j < i < k$).  Given an $H$-protrusion ${\L} = [j,k]$, the associated subgallery of $\gamma$ is
		\[
		\gamma_{{\L}} := (p_j , c_j , p_{j+1} , \ldots , p_{k-1} , c_{k-1} , p_k).
		\]
		Similarly, if $p_j\subset H$ and $c_i\not\subset H^\id$ for all $i \geq j$, then ${\L}=[j,\infty)$ is an $H$-protrusion, and we define $\gamma_{\L}$ to be the subgallery
		\[
		\gamma_{{\L}} := (p_j , c_j , p_{j+1} , \ldots , p_l , c_l ,  p_{l+1}).
		\]
	\item Call ${\L}={\L}_1 \sqcup {\L}_2 \sqcup \dots \sqcup {\L}_r \subset \N$ an \emph{$H$-outcrop (of $\gamma$)} if each ${\L}_i$ is an $H$-protrusion. 
		If ${\L}$ contains every possible $H$-protrusion, then ${\L}$ is the (unique) \emph{maximal $H$-outcrop}.  If $[j,\infty)$ is an $H$-protrusion, where $j \leq l$ is the largest index such that $p_j \subset H$, and ${\L}$ contains every possible $H$-protrusion except for $[j,\infty)$, then we call ${\L}$ the (unique) \emph{near-maximal $H$-outcrop}.  We take the maximal or near-maximal $H$-outcrop to be ${\L} = \emptyset$ if there are no $H$-protrusions.  
	\item For $1 \leq j < k \leq l+1$, call ${\L}=[j,k] \subset \N$ an \emph{$H$-indentation} if $p_j,p_k\subset H$ and $c_i \subset H^\id$ for all $j\leq i <k$; note that if $j$ is the smallest index such that $p_j \subset H$ then the interval $[0,j]$ is by definition not an $H$-indentation.
	\item Call ${\L}={\L}_1 \sqcup {\L}_2 \sqcup \dots \sqcup {\L}_r$ where each ${\L}_i$ is an $H$-indentation an \emph{$H$-ingrowth}. 
	If ${\L}$ contains every possible $H$-indentation, then ${\L}$ is the (unique) \emph{maximal $H$-ingrowth}.  We take the maximal $H$-ingrowth to be ${\L} = \emptyset$ if there are no $H$-indentations.
\end{enumerate}
\end{definition}

Note that every gallery has a maximal outcrop, but that a near-maximal outcrop does not exist if the final simplex of $\gamma$ is contained in $H^\id$. 

\begin{example}\label{Ex:OutcropsIngrowths}
Figure~\ref{fig:OutcropsIngrowths} shows an example of a gallery $\gamma = (p_0,c_1,\dots,p_{29},c_{29},p_{30})$, where $p_0 = v_0$ is the origin and $c_0 = \fa$. This gallery is illustrated with a black line. The red hyperplane $H$ contains the panels $p_5$, $p_{13}$, $p_{18}$, and $p_{26}$ of $\gamma$.  The intervals $[5,13]$, $[18,26]$, $[18,30]$, $[26, \infty]$, and $[18,\infty)$ are the $H$-protrusions, the maximal $H$-outcrop is $[5,13] \sqcup [18,\infty)$, and the near-maximal $H$-outcrop is $[5,13] \sqcup [18,26]$. The $H$-indentation $[13,18]$ is equal to the maximal $H$-ingrowth.
\end{example}

We will now define the first of the two operators in question. 

\begin{definition}\label{def:e-operators}
Let $H$ and $\gamma$ be as in Definition \ref{def:OutcropsIngrowths}. 	
Fix an $H$-protrusion ${\L}=[j,k]$. 
If ${\L} = \emptyset$, define $e^{\L}_H(\gamma) = \gamma$.
Otherwise, define $e^{\L}_H$ to be the operator which folds the subgallery $\gamma_{\L}$ onto the half-apartment $H^\id$, that is,
	\[
	e^{\L}_H(\gamma)=(p_0 , c_0 ,  \ldots , c_{j-1} , p_j , r_H(c_j) , \ldots ,  r_H(c_{k-1}) , p_k , \dots , c_l ,  p_{l+1}).
	\]
Similarly, if ${\L}=[j,\infty)$, define $e^{\L}_H$ by
	\[
	e^{\L}_H(\gamma)=(p_0 , c_0 , p_1  , \ldots , c_{j-1} , p_j , r_H(c_j) , \ldots , r_H(p_l) , r_H(c_l) ,  r_H(p_{l+1})).
	\]
For an $H$-outcrop ${\L}={\L}_1 \sqcup {\L}_2 \sqcup \ldots \sqcup {\L}_r$, we define $e^{{\L}}_H=e^{{\L}_1}_H\circ e^{{\L}_2}_H \circ \dots \circ e^{{\L}_r}_H$. Note that since the intervals ${\L}_i$ are disjoint, the operators $e^{{\L}_i}_H$ pairwise commute.	
\end{definition}

The second operator,  $f^{{\L}}_H$, is the inverse of $e^{\L}_H$, that is, the operator that unfolds a gallery across $H$.

\begin{definition}
Let $H$ and $\gamma$ be as in Definition \ref{def:OutcropsIngrowths}. 	
For ${\L} = [j,k]$ an $H$-indentation, define 
\[
f^{\L}_H(\gamma)=(p_0 , c_0 ,  \ldots , c_{j-1} , p_j , r_H(c_j) , \ldots ,  r_H(c_{k-1}) , p_k , \dots , c_l ,  p_{l+1}).
\]
Similarly, if $p_j\subset H$ with $j\geq 1$ and $c_i\subset H^\id$ for all $i \geq j$, then ${\L}=[j,\infty)$ is an {$H$-indentation}, and we define
\[
f^{\L}_H(\gamma)=(p_0 , c_0 ,  \ldots , c_{j-1} , p_j , r_H(c_j) , \ldots , r_H(p_l)  , r_H(c_l) ,  r_H(p_{l+1})).
\]

Let now ${\L}={\L}_1 \sqcup {\L}_2 \sqcup \ldots \sqcup {\L}_r$ be an $H$-ingrowth. As with the $e^{\L}_H$ operators, we define $f^{\L}_H(\gamma) = \gamma$ if ${\L} = \emptyset$ and otherwise we define $f^{\L}_H = f^{{\L}_1}_H\circ f^{{\L}_2}_H \circ \dots \circ f^{{\L}_r}_H$.  
\end{definition}

Note that by definition $f^{\L}_H$ preserves the initial simplex of $\gamma$.

\begin{remark} The operators $e^{\L}_H$ and $f_H^{\L}$ are inspired by the root operators defined in~\cite{GaussentLittelmann}, and this motivates our notation.  However, we allow folds in hyperplanes orthogonal to non-simple root directions, and $e^{\L}_H$ and $f_H^{\L}$ may reflect more than one subgallery of $\gamma$ (in the case that ${\L}$ is an $H$-outcrop or $H$-ingrowth).
\end{remark}

\begin{remark}  Given a gallery $\gamma$, the galleries $e_H^{\L}(\gamma)$ and $f_H^{\L}(\gamma)$ have the same type as $\gamma$. This is because both $e_H^{\L}$ and $f_H^{\L}$ fix the initial simplex of $\gamma$, and the reflection $r_H$ is type-preserving.
\end{remark}

The next lemma, concerning the effect of certain $e_H^{\L}$ and $f_H^{\L}$ on orientations, is the key technical result for the proof of Theorem~\ref{thm:AlcoveRecursion}.  Recall from Section~\ref{sec:chimneysRetractionsGalleries} that for $J = \{ 1,\dots, n\}$ and $y \in \aW$, the unique $(J,y)$-sector is the alcove $\y = y\fa$, and a gallery $\gamma$ is positively folded with respect to the $(J,y)$-chimney if and only if for every panel $p$ at which $\gamma$ is folded, the hyperplane supporting $p$ separates the (identical) alcoves of $\gamma$ which occur before and after this fold from $\y$.

\begin{lemma}\label{lem:GammaFold} Let $J = \{ 1,\dots,n\}$, $y \in \aW$, and $s \in S$, and let $H$ be the unique hyperplane separating the alcoves $y\fa$ and $ys\fa$.  Assume that $\ell(ys) > \ell(y)$.
\begin{enumerate}
\item If a gallery $\gamma$ is positively folded with respect to the $(J,y)$-chimney, then for the maximal or  
near-maximal $H$-outcrop ${\L}$, the gallery $e^{\L}_H(\gamma)$ is positively folded with respect to the $(J,ys)$-chimney. 
\item If a gallery $\gamma$ is positively folded with respect to the $(J,ys)$-chimney, then for the maximal $H$-ingrowth ${\L}$, the gallery $f^{\L}_H(\gamma)$ is positively folded with respect to the $(J,y)$-chimney. 
\end{enumerate}
\end{lemma}
\begin{proof} We will prove (1); the proof of (2) is similar.  Since $\ell(ys) > \ell(y)$, we observe that $H^y=H^\id$.  
Let $\gamma=(p_0 , c_0 , p_1 , c_1 , p_2 , \ldots , p_l , c_l ,  p_{l+1})$ and let ${\L}$ be the maximal $H$-outcrop.  Then 
\[
e^{{\L}}_H(\gamma) = (q_0 , d_0 , q_1 , d_1 , q_2 , \ldots , q_l , d_l ,  q_{l+1}),
\]
where $q_i=\begin{cases} p_i & \mbox{ if } p_i\subset H^y \\ r_H(p_i) & \mbox{ if } p_i\subset H^{ys}\ \end{cases} \qquad$ and $\qquad d_i=\begin{cases} c_i & \mbox{ if } c_i \subset H^y \\ r_H(c_i) & \mbox{ if } c_i \subset H^{ys}.\end{cases}$

Suppose that $e^{{\L}}_H(\gamma)$ has a fold at its panel $q_i$, so that $d_{i-1}=d_i$.  We wish to show that the hyperplane supporting the panel $q_i$ separates $d_{i-1} = d_i$ from the alcove $ys\fa$.

\emph{Case 1:} Suppose $c_{i-1}\neq c_i$, that is, $\gamma$ is not folded at $p_i$. Then $e^{\L}_H$ acts as the identity on one of the alcoves $c_{i-1}$ and $c_i$ and as $r_H$ on the other. Assume first that $d_{i-1}=c_{i-1}$ and $d_{i}=r_H(c_{i})$. Then since $c_{i-1} \subset H^{y}$ and $c_{i} \subset H^{ys}$, we have $p_i \subset H$ and $q_i =p_i$. Since $d_i=d_{i-1}=c_{i-1} \subset H^{y}$, this implies that $H$ is the hyperplane supporting $q_i$, and that $H$ separates the alcoves $d_i$ and $ys\fa$, as required.  The proof when $d_{i}=c_{i}$ and $d_{i-1}=r_H(c_{i-1})$ is similar.

\emph{Case 2:} Suppose $c_{i-1}=c_i$, that is, the gallery $\gamma$ has a fold at $p_i$.  Let $H_i$ be the hyperplane supporting $p_i$.  Since $\gamma$ is positively folded with respect to the $(J,y)$-chimney, we know that $H_i$ separates the alcoves $c_i$ and $y\fa$.  We consider two subcases.

\emph{Case 2(a):} The alcove $c_{i-1} = c_i$ is in $H^{ys}$. Then as ${\L}$ is the maximal $H$-outcrop, we have $d_{i-1}=d_i=r_H(c_i)$ and $q_i=r_H(p_i)$.  Since  the hyperplane $H_i$ separates the alcoves $c_i$ and $y\fa$, the hyperplane $r_H(H_i)$ separates the alcoves $r_H(c_i)$ and $r_H(y\fa)$.  Hence the hyperplane supporting $q_i$ separates the alcoves $d_i$ and $ys\fa$, as required.

\emph{Case 2(b):}  The alcove $c_{i-1} = c_i$ is in $H^{y}$. Then $d_{i-1}=d_i=c_i$ and $q_i=p_i$.  Since $H_i$ separates the alcoves $c_i$ and $y\fa$, it follows that $d_i \not\subset H_i^y$.   Therefore $H_i \neq H$, and so $H_i$ does not separate the alcoves $ys\fa$ and $y\fa$ (since $H$ is the unique hyperplane to do so).  Hence $ys\fa \subset H_i^y$. Thus  $H_i$, the hyperplane supporting $q_i$, separates the alcoves $d_i$ and $ys\fa$. 

This completes the proof of (1) for ${\L}$ the maximal $H$-outcrop.  Now let $j \leq l$ be the largest index such that $p_j \subset H$ and assume that $[j,\infty)$ is an $H$-protrusion, so that if ${\L}' = {\L} \setminus [j,\infty)$ then ${\L}'$ is the near-maximal $H$-outcrop.  Continuing notation from above, we have that
\[
e^{{\L}'}_H(\gamma) = (q_0 , d_0 , q_1 , d_1 , q_2 , \ldots , d_{j-1} , p_j , c_j , \dots , c_l , p_{l+1}).
\]
That is, the gallery $e_H^{{\L}'}(\gamma)$ is the same as $e_H^{\L}(\gamma)$ except that we do not apply the reflection $r_H$ to the tail $\gamma_{[j,\infty)} =  (p_j , c_j , \ldots , c_l , p_{l+1})$.  

We have already proved that $e_H^{\L}(\gamma)$ is positively folded with respect to the $(J,ys)$-chimney.  Hence the initial subgallery 
\[(q_0 , d_0 , q_1 , d_1 , q_2 , \ldots , d_{j-1} , p_j)\]
of $e_H^{{\L}'}(\gamma)$ is positively folded with respect to the $(J,ys)$-chimney.  Next, consider the consecutive alcoves $d_{j-1}$ and $c_j$.  Since $p_j$ is the last panel of $\gamma$ to lie in $H$ and $[j,\infty)$ is an $H$-protrusion, we have that $c_j \subset H^{ys}$.  As ${\L}$ is the maximal $H$-outcrop, if $c_{j-1} \subset H^{ys}$ then $d_{j-1} = r_H(c_{j-1}) \subset H^{y}$.  Also if $c_{j-1} \subset H^{y}$ then $d_{j-1} = c_{j-1} \subset H^{y}$.  Hence $d_{j-1} \neq c_j$ and so the gallery $e^{{\L}'}_H(\gamma)$ is not folded at $p_j$.  
Since $p_j$ is the last panel of $\gamma$ to be contained in $H$, any folds in the tail $\gamma_{[j,\infty)}$ are in hyperplanes $H' \neq H$.  Now as the alcoves $ys\fa$ and $y\fa$ are not separated by any hyperplanes except for $H$, it follows that $\gamma_{[j,\infty)}$ is positively folded with respect to the $(J,ys)$-chimney.  We conclude that $e_H^{{\L}'}(\gamma)$ is positively folded with respect to the $(J,ys)$-chimney, as required.
\end{proof}

Note that in item (1) of Lemma~\ref{lem:GammaFold} the gallery $e^{\L}_H(\gamma)$ is positively folded with respect to both the maximal and near-maximal outcrop in case both are defined. 

We are now ready to prove the main result of this section.  Although Theorem~\ref{thm:AlcoveRecursion} below is stated only for shadows with respect to $(J,y)$-chimneys where $J = \{1,\dots,n\}$ (see Notation~\ref{not:shadow}), by combining Corollary~\ref{cor:sim-stable-alcove} and Proposition~\ref{prop:shadowsRetractions} with Theorem~\ref{thm:AlcoveRecursion} we obtain a recursive description of shadows for all chimneys, as the sets $r^{-1}(W_\sigma \cdot x\tau)$, $r^{-1}(\x)$, and $r^{-1}(\sW \cdot \lambda)$ appearing in Proposition~\ref{prop:shadowsRetractions} are all bounded.  The key point in the proof of Theorem~\ref{thm:AlcoveRecursion} is that galleries which are positively folded with respect to the $(J,y)$-chimney may not be positively folded with respect to the $(J,ys)$-chimney, and so some careful manipulation of galleries is needed in order to calculate the shadow correctly.

\begin{thm}\label{thm:AlcoveRecursion}  Let $x \in \aW$ and let $\sigma$ and $\tau$ be faces of the fundamental alcove $\fa$ which contain the origin.  Let $y \in \aW$ and $s \in \aS$.  If $\ell(ys)>\ell(y)$, then 
\[
\Sh_{\ys}(x\tau,\sigma)=\Sh_{\y}(x\tau,\sigma)\cup r_H ( H^{ys} \cap \Sh_\y(x\tau,\sigma)),
\]
where $H$ is the unique hyperplane separating the alcoves $ys\fa$ and $y\fa$.
\end{thm}

\begin{proof} As in the proof of Lemma~\ref{lem:GammaFold}, the hypothesis $\ell(ys) > \ell(y)$ implies that $H^{y} = H^\id$.  Fix a minimal gallery $\gamma_x$ with first face $\sigma$ and last face $x\tau$.  We will use $\Gamma^+_{\y}(\gamma_x,\zeta)$ to denote the set of all galleries of the same type as $\gamma_x$ which end at the simplex $\zeta$ and are positively folded with respect to the $(J,y)$-chimney, where $J = \{1,\dots,n\}$.

Suppose first that we have $\zeta \in \Sh_{\y}(x\tau,\sigma)$.  Then by definition of the shadow, there is a gallery $\gamma \in \Gamma^{+}_{\y}(\gamma_x,\zeta)$.  Let ${\L}$ be the maximal $H$-outcrop of $\gamma$ if $\zeta \in H^{y}$ and the near-maximal $H$-outcrop of $\gamma$ if $\zeta \in H^{ys}$ (if $\zeta \in H$, then either choice will be valid).  Then by Lemma~\ref{lem:GammaFold}(1), the gallery $e_H^{\L}(\gamma)$ is positively folded with respect to the $(J,ys)$-chimney.  For ${\L}$ both maximal and near-maximal, by choice of ${\L}$ the final simplex of $e_H^{\L}(\gamma)$ is $\zeta$.   Hence $e_H^{\L}(\gamma) \in \Gamma^{+}_{\ys}(\gamma_x,\zeta)$, and so $\zeta \in \Sh_{\ys}(x\tau,\sigma)$.  

Now suppose we have $\eta \in r_H(H^{ys} \cap \Sh_{\y}(x\tau,\sigma))$.  Then there exists $\gamma \in \Gamma^{+}_{\y}(\gamma_x,\zeta)$, where $\zeta = r_H(\eta)$, equivalently $\eta = r_H(\zeta)$, and $\zeta \in H^{ys}$.  Since the first simplex of $\gamma$ is $\sigma \subset \fa$ (by definition of the shadow), it follows that there is at least one panel of $\gamma$ contained in $H$.  So letting $\gamma$ have final simplex $p_{n+1} = \zeta$, the maximal $H$-outcrop ${\L}$ includes the index $n+1$.  Thus $e_H^{\L}(\gamma)$ has final simplex $r_H(\zeta) = \eta$.  By Lemma \ref{lem:GammaFold}(1), we have that $e_H^{\L}(\gamma)$ is positively folded with respect to the $(J,ys)$-chimney.  So $e_H^{\L}(\gamma) \in \Gamma^{+}_{\ys}(\gamma_x,\eta)$.  Hence $\eta \in \Sh_{\ys}(x\tau,\sigma)$.

For the other inclusion, let $\zeta$ be in $\Sh_{\ys}(x\tau,\sigma)$ but not in $\Sh_{\y}(x\tau,\sigma)$. Then there exists a gallery $\gamma$ from $\sigma$ to $\zeta$ that is positively folded with respect to the $(J,ys)$-chimney, but is \emph{not} positively folded with respect to the $(J,y)$-chimney.  Let ${\L}$ be the maximal $H$-ingrowth of $\gamma$.  By Lemma~\ref{lem:GammaFold}(2), the gallery $f_H^{\L}(\gamma)$ is positively folded with respect to the $(J,y)$-chimney.  If $f_H^{\L}(\gamma)$ ends in $\zeta$, then we contradict the fact that $\zeta \notin \Sh_{\y}(x\tau,\sigma)$. So $f_H^{\L}(\gamma)$ ends in $r_H(\zeta)$, implying that $\zeta \in H^{y}$, and hence $r_H(\zeta) \in H^{ys} \cap \Sh_{\y}(x\tau,\sigma)$.  Thus $\zeta \in r_H(H^{ys} \cap \Sh_{\ys}(x\tau,\sigma))$, which completes the proof.
\end{proof}

\begin{remark} \label{rem:braidRecursion} We observe that the proof of Theorem~\ref{thm:AlcoveRecursion} shows that for any choice of minimal gallery $\gamma_x$ from $\fa$ to $\x$, and for every alcove $\z$ in the shadow $\Sh_{\ys}(\x)$, there exists a gallery in $\Gamma^+_{\y}(\gamma_x,\z)$. In particular, given two different braid-equivalent reduced words for $x$, we have distinct minimal galleries $\gamma_x$ and $\gamma'_x$ with types given by these reduced words. Thus, for any $z\in W$ such that $\z \in \Sh_{\ys}(\x)$, we obtain galleries in both $\Gamma^+_{\y}(\gamma_x,\z)$ and $\Gamma^+_{\y}(\gamma'_x,\z)$.  This gives us another proof that the orientation induced by the $(J,y)$-chimney for $J=\{1,\dots,n\}$ is braid-invariant.  Then by applying Corollary~\ref{cor:sim-stable-alcove}, we obtain another proof of braid-invariance for general chimneys (see Proposition~\ref{prop:braidinvariant}).
\end{remark}

\begin{remark}\label{rem:multiplicities}  The proof of Theorem~\ref{thm:AlcoveRecursion} calculates the shadow using sets of positively folded galleries of the form $\Gamma^+_{\y}(\gamma_x,\zeta)$.  It is reasonable to ask whether one can determine the number of distinct galleries in $\Gamma^+_{\y}(\gamma_x,\zeta)$, that is, whether one can calculate a  ``shadow with multiplicities".  Determining such multiplicities is a delicate problem, essentially because the operators $e_H^L$ and $f_H^L$ used in the proof of Theorem~\ref{thm:AlcoveRecursion} may reflect all or only some of the portions of a gallery lying on the appropriate side of $H$.
\end{remark}

As a corollary of Theorem~\ref{thm:AlcoveRecursion}, we obtain the following restriction on where the shadow of a coroot lattice element $\lambda$ can lie.  In the case that $X$ is the Bruhat--Tits building for $G(F)$ with standard apartment $A$,  and the retraction onto $A$ centered at $\y$ has the same effect on $r^{-1}(\sW \cdot \lambda)$ as the retraction from the chamber at infinity corresponding to the negative unipotent subgroup $U^-$ of $G(F)$, Corollary~\ref{cor:WeylPolytope} recovers a well-known containment found in Bruhat--Tits \cite{BruhatTits}. The convex hull $\Conv$ in the following statement is the metric convex hull in the apartment $A$.

\begin{corollary}\label{cor:WeylPolytope}  Let $\lambda \in R^\vee$ and $y \in \aW$.  Then 
\[
\Sh_\y(\lambda) \subseteq \Conv(\sW \cdot \lambda).
\]
That is, for $J = \{1,\dots,n\}$ and $y \in \aW$, the shadow of $\lambda$ with respect to the $(J,y)$-chimney is contained in the $\lambda$-Weyl polytope.
\end{corollary}

\begin{proof}  The proof is by induction on $\ell(y)$.  If $\ell(y) = 0$, that is, $y = \id$, then $\Sh_\y(\lambda)$ equals the orbit $\sW \cdot \lambda$ by Proposition \ref{prop:shadowsRetractions}, and we are done.  Now let $s_{i_1} \dots s_{i_k}$ be a reduced word for $y$.  To simplify notation, let $y' = s_{i_1} \dots s_{i_{k-1}}$ and $s = s_{i_{k}}$, and let $H$ be the hyperplane separating the alcoves $y'\fa$ and $y's\fa = y\fa$.  By Theorem~\ref{thm:AlcoveRecursion} and inductive assumption, it suffices to show that the set $r_H(H^{y's} \cap \Conv(\sW \cdot \lambda))$ is contained in $\Conv(\sW \cdot \lambda)$.  If the intersection $H^{y's} \cap \Conv(\sW \cdot \lambda)$ is empty there is nothing to show, so we assume that there is some $\mu \in H^{y's} \cap \Conv(\sW \cdot \lambda)$.  

As the reflection $r_H$ fixes the hyperplane $H$, by considering convex hulls it is now enough to show that for each $\mu \in H^{y's} \cap (\sW \cdot \lambda)$, we have $r_H(\mu) \in \Conv(\sW \cdot \lambda)$.  Moreover, by applying an appropriate element of $\sW$, it suffices to take $\mu = \lambda \in H^{y's}$.  That is, we are in the situation that the hyperplane $H$ separates the origin from $\lambda$, and we wish to show that $r_H(\lambda)$ lies in $\Conv(\sW \cdot \lambda)$.  

Now $H = H_{\beta,i}$ for some root $\beta$ (not necessarily simple) and some index $i$, where without loss of generality $i > 0$.  More precisely, letting $k_\lambda = \langle \beta, \lambda \rangle$ we have without loss of generality that $0 < i < k_\lambda$, since the origin lies in $H_\beta = H_{\beta,0}$, the coroot lattice point $\lambda$ lies in the hyperplane $H_{\beta,k_\lambda}$, and $H = H_{\beta,i}$ separates the origin from $\lambda$.  The reflection $s_\beta$ in the hyperplane $H_{\beta}$ takes $\lambda$ to another extremal vertex of $\Conv(\sW \cdot \lambda)$, namely the vertex $s_{\beta}(\lambda) = \lambda - k_\lambda \beta^\vee$.  The image of $\lambda$ under the reflection $s_{\beta,i}$ in $H_{\beta,i}$ then lies on the line interval connecting $\lambda$ and $s_\beta(\lambda)$, since $s_{\beta,i}(\lambda) = \lambda - (k_\lambda - i)\beta^\vee$ and $0 < i < k_\lambda$.  This line interval is contained in $\Conv(\sW \cdot \lambda)$ and we have $r_H = s_{\beta,i}$, so $r_H(\lambda)$ is contained in $\Conv(\sW \cdot \lambda)$ as required.
\end{proof}

In the following lemma we prove a more general version of the assertion of Corollary~\ref{cor:WeylPolytope}. As we will not use it in its full generality we will not formally introduce some of the notions used in the proof, but refer the reader to \cite{GraeberSchwer} instead. 

\begin{lemma}
Let $\lambda \in R^\vee$ and let $\phi$ be any orientation (not just one induced by a chimney). Then 
$	\Sh_\phi(\lambda) \subseteq \Conv(\sW \cdot \lambda).$
That is, the shadow of a coroot lattice element $\lambda$ with respect to any orientation is contained in the $\lambda$-Weyl polytope.
\end{lemma}
\begin{proof}
	The shadow with respect to the trivial positive orientation contains (by definition) the end-vertices of all folded galleries of a fixed type. One can show that these end-vertices are just those elements $\nu$ of the coroot lattice whose dominant image $\nu^+$ under the $\sW$-action is smaller than $\lambda^+$, i.e. all $\nu$ such that $\nu^+ \leq \lambda^+$ in dominance order. This follows from the fact that Bruhat order on group elements (corresponding to all the end-alcoves of folded galleries of a fixed type) restricts to dominance order on vertices.  But the condition  $\nu^+ \leq \lambda^+$ in dominance order is the same as taking the convex hull of $\sW\cdot\lambda$ as shown in, for instance, Lemma 3.5 of~\cite{HitzelThesis}. Any vertex shadow is a subset of the trivial positive shadow and hence a subset of the $\lambda$-Weyl polytope.   
\end{proof}


\section{Positively folded galleries and double cosets}\label{sec:galleriesCosets}

In this section we establish the connection between the previous combinatorial results on folded galleries, and double coset intersections inside the underlying algebraic groups in the function field case.  We first record some additional definitions and notation in Section~\ref{sec:conventions}, then in Section~\ref{sec:Pchimney} adapt the key concepts from Section~\ref{sec:chimneysRetractionsGalleries} to this setting.  Section~\ref{sec:DoubleCosets} proves the main result of this section, Theorem~\ref{thm:DoubleCosetIntersectionParahoric}, which generalizes Theorems~\ref{thm:DoubleCosetIntersectionI} and~\ref{thm:Grassmannian-intro} of the introduction.  This theorem establishes a bijection between the points of certain double coset intersections in partial affine flag varieties, and certain sets of positively folded galleries which have been further ``decorated" by elements of the residue field.  We then prove a nonemptiness statement in Theorem~\ref{thm:IntersectionParahoricShadow}, and finally apply Theorems~\ref{thm:DoubleCosetIntersectionParahoric} and~\ref{thm:IntersectionParahoricShadow} to the case of the affine Grassmannian in Corollary~\ref{cor:DoubleCosetIntersectionK}.

\subsection{Conventions}\label{sec:conventions}

We continue all terminology and notation from Section~\ref{sec:preliminaries}, and include here some further items needed to state and prove the theorems in this section.  We are not redefining anything from the Coxeter-theoretic setup used in Section~\ref{sec:preliminaries}, but rather providing additional algebraic interpretations for  many of these notions.

Let $\res$ be any field (not necessarily finite), let $F = \laurent$, and let $G$ be a split, connected, reductive group over $F$. Let $\mathcal{O} = \res[[t]]$ be the ring of integers of $F$.  Choose a split maximal torus $T$ of $G$ and a Borel subgroup $B = TU$, where $U$ is the unipotent radical. Denote the underlying root system by $\Phi$, and further assume that $\Phi$ is irreducible.  Throughout the remainder of the paper, we denote by $\Phi^+$ those (spherical) roots which are positive with respect to the \emph{opposite} Borel subgroup $B^-$.  We will discuss this convention further below.

Let $\sW$ be the (finite or spherical) Weyl group of $T$ in $G$; that is, $\sW = N_G(T) / T$. For $\lambda \in X_*(T)$, write $t^{\lambda}$ for the element in $T(F)$ that is the image of $t$ under the homomorphism $\lambda: \mathbb{G}_m \rightarrow T$.   Let $\Delta = \{\alpha_i \}_{i=1}^n$ be a basis of simple roots in $X^*(T)$, chosen such that each $\alpha_i \in \Phi^+$ is positive with respect to $B^-$. We denote by $\alpha^{\vee} = 2 \alpha / \langle \alpha, \alpha \rangle$ the coroot associated to $\alpha \in \Phi$ with respect to the pairing $\langle \ ,\  \rangle: X^*(T) \times X_*(T) \rightarrow \mathbb{Z}$. 
Let  $R^{\vee} = \oplus_{i=1}^n \Z \alpha_i^{\vee} \subset X_*(T)$ be the coroot lattice.  
 Let $\aW = R^{\vee} \rtimes \sW$ be  the affine Weyl group, and $\widetilde{W} = X_*(T) \rtimes \sW$ be the extended affine Weyl group.  
We view $\sW$  and $\aW$ as acting on the vector space $V = \R^n$, which we can identify with $X_*(T)\otimes_{\Z} \R$.
By abuse of notation, we use the same symbol for an element of $\sW$ or $\aW$ and a lift of this element to $G(F)$.

Let $X$ be the Bruhat--Tits building for $G(F)$.  Then $X$ is a (thick) affine building of type $(W,S)$.  Our choice of split maximal torus $T$ determines the standard apartment $\App$ of $X$.  Under the identification $\App = X_*(T) {\otimes}_{\Z} \R = V$, the interior of any alcove in $\App$ can be described as the set of all points $v \in V$ satisfying the strict inequalities $k_{\alpha} < \langle \alpha, v \rangle < k_{\alpha} + 1$, where $\alpha$ runs through~$\Phi^+$ and $k_{\alpha} \in \Z$ is some fixed integer (depending upon the alcove and $\alpha$).  Recall that the dominant Weyl chamber $\Cf$ is the set of points $v \in V$ such that $\langle \alpha, v \rangle \geq 0$ for all $\alpha \in \Phi^+$.  We note that our convention of $\Phi^+$ being the (spherical) roots which are positive with respect to $B^-$ in fact ensures that the elements of $\Phi$ which lie in the half-apartment of $\App$ which is bounded by the hyperplane $H_{\widetilde\alpha}$ and contains the dominant Weyl chamber are exactly the positive roots $\Phi^+$.

As in~\cite{GHKRadlvs}, we choose the Iwahori subgroup $I$ of $G(F)$ to be the inverse image of the opposite Borel subgroup $B^-(\res)$ under the projection map $G(\mathcal{O}) \rightarrow G(\res)$.  This convention ensures that for all $\lambda$ in the coroot lattice, the element $t^\lambda$ acts on $\App$ by translation by $\lambda$ (rather than by $-\lambda$, had we chosen $I$ to be the preimage of $B$).  The base alcove $\fa$ of $\App$ is the (unique) alcove of $X$ whose stabilizer in $G(F)$ is equal to~$I$.

Given any $\alpha \in \Phi$ and $k \in \Z$, we write $U_{\alpha,k}$ for the subgroup of $G(F)$ which fixes pointwise the half-apartment of $\App$ which is bounded by the hyperplane $H_{\alpha,k}$ and contains the points $v\in A$ such that $\langle \alpha, v \rangle \geq k$.  (Thus if $\alpha \in \Phi^+$ then $U_{\alpha,k}$ fixes the half-apartment $\alpha^{k,\id}$ defined in Section~\ref{sec:buildings}, while if $\alpha \in \Phi^-$ then $U_{\alpha,k}$ fixes the half-apartment $\alpha^{k,w_0}$.)
The pair $(\alpha, k)$ corresponds to a unique affine root $\beta$ and vice versa, so we sometimes denote the subgroup $U_{\alpha,k}$ by $U_\beta$, and the reflection in the hyperplane $H_{\beta}:=H_{\alpha,k}$ by $s_\beta \in \aW$.  Recall that $\aW$ acts on the set of affine roots.  Given $x \in \aW$ and an affine root~$\beta$, we write $x(\beta)$ for the image of $\beta$ under $x$.  We also write $\alpha_0$ for the affine root such that $U_{\alpha_0} = U_{-\widetilde\alpha,-1}$, so that $U_{\alpha_0}$ fixes the half-apartment consisting of the points $v \in A$ such that $\langle -\widetilde\alpha,v\rangle \geq -1$, or equivalently $\langle \widetilde\alpha,v\rangle \leq 1$.  In particular, $U_{\alpha_0}$ fixes the base alcove and thus $U_{\alpha_0}$ is a subgroup of $I$.  We write $s_0$ for the affine reflection~$s_{\widetilde\alpha,1}$.

\subsection{$P$-chimneys and the associated retractions}\label{sec:Pchimney}

Let $P=P(F)$ be a standard spherical parabolic subgroup of $G(F)$, meaning a subgroup of $G(F)$ which contains the fixed Borel subgroup $B = B(F)$.  Then $P$ has Levi decomposition $P = MN$, where $M$ is the Levi component and $N$ is its unipotent radical, with both $M$ and $N$ groups over $F$.  Define $\Delta_M \subseteq \Delta$ to be the set of positive simple roots $\alpha_j$ such that the subgroup $U_{\alpha_j}$ is contained in $M$.  Recall that in our conventions, each $\alpha_j$ is positive with respect to the opposite Borel subgroup $B^-$; however, note that for any $\alpha \in \Phi$, the subgroup $U_\alpha$ is contained in $M$ if and only if $U_{-\alpha}$ is contained in $M$.  Let $\Phi_M \subseteq \Phi$ be the sub-root system generated by $\Delta_M$ and define $\Phi_M^+ = \Phi^+ \cap \Phi_M$.  Then $\Phi_M^+$ consists of the elements $\alpha \in \Phi$ which are positive with respect to $B^-$ and are such that $U_\alpha$ is contained in $M$.  We remark that putting $\alpha' = -\alpha$, this is equivalent to $\alpha' \in \Phi$ being positive with respect to $B = B^+$ and $U_{\alpha'}$ being contained in $M$.
    
Now let $J \subseteq \{ 1,\dots,n\}$ be such that $\Delta_M = \{ \alpha_j \mid j \in J \}$.  Then
\[ P = P_J = \bigsqcup_{w \in W_J} BwB,\]
where $W_J$ is the standard parabolic subgroup of $\sW$ generated by $\{ s_j \mid j \in J \}$.
Finally, the subgroup $I_P$ of $G(F)$ is defined to be
\[
I_P := (I \cap M)N = (I \cap M(F))N(F).
\]
Note that since $I$ is a subgroup of $G(\cO)$, the intersection $I \cap M$ will be a subgroup of $G(\cO)$ as well, while $N$ is defined over  $F$.  Also, the spherical subgroups $U_{\alpha,0}$ which are contained in $N$ will be all those of the form $U_{-\alpha}$ where $\alpha \in \Phi \setminus \Phi_M$ is positive with respect to $B^-$.  In other words, as one would expect, those subgroups which are contained in $N$ will be all those  of the form $U_{\alpha'}$ where $\alpha' \in \Phi \setminus \Phi_M$ is positive with respect to $B$. We denote by $R(I_P)$ the set of all affine roots $\beta$ such that the subgroup $U_\beta$ is contained in $I_P$.  Note that $\beta \in R(I_P)$ if and only if $-\beta \notin R(I_P)$.

The following Bruhat decomposition result is established in~\cite{GHKRadlvs} in the case $F = \overline{\F_q}((t))$, where $\overline{\F_q}$ is an algebraic closure of the finite field~$\F_q$, but the proof given there extends immediately to our setting $F = \laurent$.  Recall that we denote the extended affine Weyl group of $G(F)$ by $\eW$.

\begin{lemma}[Lemma~11.2.1 of~\cite{GHKRadlvs}] Let $P = P(F)$ be a standard spherical parabolic subgroup of $G(F)$.  Then 
\[ G(F) = \bigsqcup_{z \in \widetilde{W}} I_P z I. \] 
\end{lemma}

We will need the following corollary in our proofs in Section~\ref{sec:DoubleCosets}.

\begin{corollary}\label{cor:disjoint}  Let $P = P(F)$ be a standard spherical parabolic subgroup of $G(F)$.   If $z \neq z'$ are distinct elements of the affine Weyl group $\aW$, then the double cosets $I_P z I$ and $I_P z' I$ are disjoint.
\end{corollary}

\begin{proof} This holds since $\aW$ is a subgroup of $\widetilde{W}$.
\end{proof}

We now adapt our main definitions from Section~\ref{sec:chimneysRetractionsGalleries} to this setting.  

\begin{definition}\label{def:P-chimney}  
Fix a standard spherical parabolic subgroup $P= P_J$ and let $y \in \aW$.  
The \emph{$P$-chimney}, denoted $\xi_P$, is the $J$-chimney $\xi_J$ in the standard apartment from Definition~\ref{def:J-chimney}.
For any collection $\{n_\beta\in\Z \mid \beta\in \Phi^+\setminus \Phi_M^+\}$, the corresponding \emph{$P$-sector} is $S_P(\{n_\beta\}) = S_J(\{n_\beta\})$.  
Write $S_{P}(0)$ for the $P$-sector defined by the all-zeros sequence; i.e.~$n_\beta=0$ for all $\beta\in \Phi^+\setminus \Phi_M^+$.  The \emph{$(P,y)$-chimney} is the $(J,y)$-chimney $\xi_{P,y} = \xi_{J,y}$ from Definition~\ref{def:Jy-chimney}; that is, $\xi_{P,y} = y \cdot \xi_P$.  For any collection $\{n_\beta\in\Z \mid \beta\in \Phi^+\setminus \Phi_M^+\}$, the $(P,y)$-sector $S_{P,y}(\{n_\beta\})$ is the $(J,y)$-sector $S_{J,y}(\{n_\beta\})$.  

The \emph{$(P,y)$-chimney retraction} is the retraction $r_{P,y} = r_{J,y}$ from Definition~\ref{def:chimney-retraction}.  We may put $r_P = r_J$ for the $P$-chimney retraction.  
The \emph{orientation induced by the $(P,y)$-chimney}, denoted $\phi_{P,y}$, is the orientation $\phi_{J,y}$ from Definition~\ref{def:chimneyOrientation}, and a gallery is \emph{positively folded with respect to the $(P,y)$-chimney} if it is positively folded with respect to the orientation $\phi_{P,y}$.  If $y = \id$, put $\phi_{P,y} = \phi_P$.

Let $x \in \aW$ and let $\sigma$ and $\tau$ be faces of the fundamental alcove $\fa$ which contain the origin.  
The \emph{shadow of $x\tau$ starting at $\sigma$ with respect to the $(P,y)$-orientation}, denoted $\Sh_{P,y}(x\tau, \sigma)$, is the shadow $\Sh_{J,y}(x\tau,\sigma)$ from Definition~\ref{def:shadow}.  For $\lambda \in R^\vee$, the \emph{shadow of $\lambda$ with respect to the $(P,y)$-orientation}, denoted $\Sh_{P,y}(\lambda)$, is the shadow $\Sh_{J,y}(\lambda)$ from Notation~\ref{not:shadow}.
\end{definition}

\begin{remark}\label{rem:GHKR-retraction}  The retraction $r_{P,y}$ generalizes the retraction $\rho_{P,w}$ described in Section 11.2 of~\cite{GHKRadlvs} in the case that $F = \overline{\F_q}((t))$, where $\overline{\F_q}$ is an algebraic closure of the finite field~$\F_q$.  For $P = MN$, it is stated in~\cite{GHKRadlvs} that the retraction $\rho_P = \rho_{P,\id}$ has the same effect as retracting from any alcove which lies between $H_\alpha$ and $H_{\alpha,1}$ for all $\alpha \in \Phi_M^+$ and is ``sufficiently  antidominant" for all roots in $N$.  We have formalized this idea using $P$-sectors.  The subgroup $I_P$ of $G(F)$ is defined exactly as above in~\cite{GHKRadlvs}.

The retraction $\rho_P$ from~\cite{GHKRadlvs} was also considered by Haines, Kapovich, and Millson in Section 6.3 of~\cite{HKM}, where this retraction is denoted $\rho_{I_P,\mathcal{A}}$.  Here, for $\Delta_M$ the $M$-dominant Weyl chamber, they composed $\rho_P:X \to \App$ with a further retraction $\App \to \Delta_M$ to obtain a retraction of the building $X$ onto $\Delta_M$.  Our $P$-sectors are essentially the same as the sets of alcoves satisfying (a) and (b) in the statement of Lemma~6.4 of~\cite{HKM}, and Lemma~6.4 of~\cite{HKM} is thus essentially the case $r_{J,y} = \rho_P$ of Corollary~\ref{cor:sim-stable-alcove} above (although our proof is different, since in Section~\ref{sec:chimneysRetractionsGalleries} we are in the setting of an arbitrary affine building).
\end{remark}

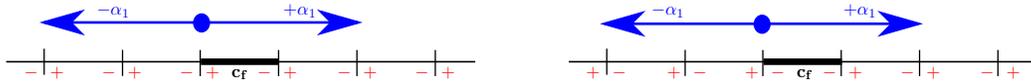
\begin{figure}[ht]
\begin{center}
\resizebox{0.9\textwidth}{!}
{
\begin{overpic}{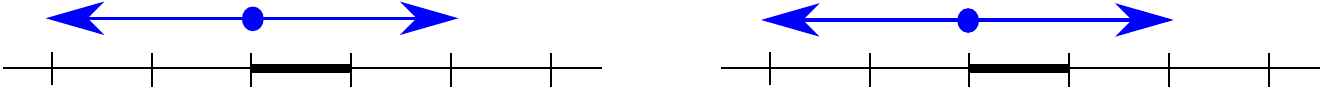}
\put(27,6){\textcolor{blue}{$+\alpha_1$}}
\put(9,6){\textcolor{blue}{$-\alpha_1$}}
\put(22,0){$\fa$}
\put(2,0){\textcolor{red}{$-$}}
\put(4.5,0){\textcolor{red}{$+$}}
\put(9.5,0){\textcolor{red}{$-$}}
\put(12,0){\textcolor{red}{$+$}}
\put(17,0){\textcolor{red}{$-$}}
\put(19.5,0){\textcolor{red}{$+$}}
\put(24.5,0){\textcolor{red}{$-$}}
\put(27,0){\textcolor{red}{$+$}}
\put(32.2,0){\textcolor{red}{$-$}}
\put(34.7,0){\textcolor{red}{$+$}}
\put(39.7,0){\textcolor{red}{$-$}}
\put(42.2,0){\textcolor{red}{$+$}}
\put(81,6){\textcolor{blue}{$+\alpha_1$}}
\put(62.5,6){\textcolor{blue}{$-\alpha_1$}}
\put(76.5,0){$\fa$}
\put(56.2,0){\textcolor{red}{$+$}}
\put(58.7,0){\textcolor{red}{$-$}}
\put(64,0){\textcolor{red}{$+$}}
\put(66.5,0){\textcolor{red}{$-$}}
\put(71.5,0){\textcolor{red}{$+$}}
\put(74,0){\textcolor{red}{$-$}}
\put(79,0){\textcolor{red}{$-$}}
\put(81.5,0){\textcolor{red}{$+$}}
\put(86.5,0){\textcolor{red}{$-$}}
\put(89,0){\textcolor{red}{$+$}}
\put(94.2,0){\textcolor{red}{$-$}}
\put(96.7,0){\textcolor{red}{$+$}}
\end{overpic}
}
\caption{Chimney-induced orientations in type $\tilde{A}_1$. For details see Example~\ref{ex:orientations}.}
\label{fig:Orientations}
\end{center}
\end{figure}

\begin{example}\label{ex:orientations}
	Two one-dimensional examples for chimney-induced orientations are provided in Figure~\ref{fig:Orientations}, where the root system $\Phi = \{ \pm \alpha_1 \}$ is of type $A_1$.  If $P = B$, equivalently $J = \emptyset$, then $\alpha_1 \in \Phi^+\setminus \Phi_M^+$ and the chimney-induced orientation is the periodic orientation on the $\alpha_1$-hyperplanes depicted on the left.  The $P$-chimney can then be identified with the point at infinity at the left-hand (antidominant) end of the apartment shown.  
	If $P = G(F)$, equivalently $J = \{ 1\}$, then $\alpha_1 \in \Phi_M^+$, and the orientation induced by the $P$-chimney is as shown on the right, so that the positive side of each panel in an $\alpha_1$-hyperplane faces away from the base alcove $\fa$; in this case, the alcove $\fa$ is the unique $P$-sector.  
\end{example}

More generally, we have the following examples of $P$-chimneys and the corresponding subgroups $I_P$. 

\begin{example}\label{eg:Bchimney}
If $P = B = B(F) = T(F)U(F)$, then $\Phi_M =\emptyset$ and we have $M = T(F)$ and $N = U(F)$, so $I_P=I_B = T(\cO)U(F)$.  As explained in Example~\ref{eg:emptyset}, the chimney $\xi_{B} = \xi_\emptyset$ can be identified with the chamber at infinity represented by the antidominant Weyl chamber.
\end{example}

\begin{example}\label{eg:GFchimney}
If $P = G(F)$ then $\Phi_M=\Phi$ and $I_P = I$.   The $P$-chimney is the collection of all half-apartments containing the base alcove $\fa$, and the only $P$-sector is the base alcove $\fa$ itself. 
\end{example}

We mention that Section~11.2 of~\cite{GHKRadlvs} includes, in effect, an argument that $I_P$ stabilizes the $P$-chimney.  In general, $I_P$ will be a proper subgroup of the stabilizer in $G(F)$ of the $P$-chimney.  For example, the stabilizer of the $B$-chimney is $B(F) = T(F)U(F)$, while $I_B = T(\cO)U(F)$.

\subsection{Double cosets}\label{sec:DoubleCosets} 

In this section, we establish several bijections between cosets in flag varieties and folded galleries, leading up to our main result, Theorem~\ref{thm:DoubleCosetIntersectionParahoric}.  Our proofs here rely on and substantially generalize the cases that were treated in \cite{PRS} and \cite{MST1}.   We consider double cosets in the affine flag variety in Proposition~\ref{prop:doublecoset} and Theorem~\ref{thm:doublecosety},  then use these results to prove Theorem~\ref{thm:DoubleCosetIntersectionParahoric}.  We next relate nonemptiness of double coset intersections to shadows, in Theorem~\ref{thm:IntersectionParahoricShadow}, and finally establish Corollary~\ref{cor:DoubleCosetIntersectionK}, which applies the results of this section to the affine Grassmannian. 

Throughout this section, $P=P(F)$ is a fixed standard spherical parabolic subgroup of $G(F)$, $y \in \aW$, and $\phi = \phi_{P,y}$ is the orientation induced by the $(P,y)$-chimney.  We now extend the definitions of a labeled folded alcove walk from~\cite{PRS} and~\cite{MST1} to all orientations induced by chimneys. 

\begin{definition} Let $\gamma = (p_0,c_0,p_1,\dots, c_l,p_{\ell+1})$ be a combinatorial gallery in the standard apartment $\App$.  For $1 \leq i \leq \ell$ write $H_i$ for the hyperplane supporting the panel $p_i$.   We say that $\gamma$ \emph{crosses $p_i$ from the negative side to the positive side (with respect to $\phi$)} if the alcove $c_{i-1}$ is on the negative side of $H_i$ with respect to $\phi$, while the alcove $c_i \neq c_{i-1}$ is on the positive side of $H_i$ with respect to $\phi$ (see Definitions~\ref{def:chimneyOrientation} and~\ref{def:positive-negative}).  The definition of $\gamma$ crossing $p_i$ from the positive side to the negative side is similar.
\end{definition}

Recall our definition of the type of a gallery (Definition~\ref{def:type}).   

\begin{definition}\label{def:labeledwalks}   Let $x \in W$ and fix a reduced word for $x$.  A \emph{labeled folded alcove walk of type $\vec{x}$ which is positively folded with respect to the $(P,y)$-chimney} is a gallery \[\gamma = (p_0,c_0,p_1,\dots, c_l,p_{\ell+1})\] 
with first face $p_0 = c_0$ and last face $p_{\ell + 1} = c_l$ 
of type $\vec{x}$ which is positively folded with respect to the $(P,y)$-chimney, together with the following labelings of the panels $p_1,\dots,p_l$:
\begin{enumerate}
\item if $\gamma$ crosses $p_i$ from the negative side to the positive side, then the panel $p_i$ is labeled by an element of $\res$; 
\item if $\gamma$ crosses $p_i$ from the positive side to the negative side, then the panel $p_i$ is labeled by $0$; and 
\item if $\gamma$ has a (positive) fold at $p_i$, then the panel $p_i$ is labeled by an element of $\res^\times$.
\end{enumerate}
The \emph{type} of a labeled folded alcove walk is the type of the underlying (unlabeled) positively folded gallery.  
\end{definition}

We note that the previous two definitions are combinatorial in nature.  Our proofs will rely upon the relationships between these combinatorial data and the groups and group elements we now define.  Recall that $U_{\alpha,k}$ denotes the subgroup of $G(F)$ which fixes pointwise the half-apartment of $\App$ which is bounded by the hyperplane $H_{\alpha,k}$ and contains the points $v\in A$ such that $\langle \alpha, v \rangle \geq k$.  Given $\alpha \in \Phi$ and $k \in \Z$, we now define the quotient \[\overline U_{\alpha,k} = U_{\alpha,k}/U_{\alpha,k+1}.\]  For $\beta$ the unique affine root corresponding to $(\alpha, k)$,  we write $\overline U_\beta = \overline U_{\alpha,k}$. For any $\beta$, the group $\overline U_\beta$ is naturally isomorphic to the additive group $(\res,+)$.  Moreover, $\overline U_\beta$ is naturally isomorphic to a subgroup of $G(F)$; we shall use this isomorphism without comment henceforth.  For example, since $\overline U_{\alpha_0}$ fixes the base alcove, we view $\overline U_{\alpha_0}$ as a subgroup of $I$. Note that throughout in the literature, the distinct groups $U_\beta$ and $\overline U_\beta$ are both frequently called root (sub)groups of $G(F)$.  To avoid confusion, we do not repeat this terminology here, and instead  use these two separate notations.

We write $u_\beta(c)$ for an element of $\overline U_\beta$, where $c \in \res$, so that $c \mapsto u_\beta(c)$ is an isomorphism from the additive group $(\res,+)$ onto $\overline U_\beta$.  
Finally, for any affine root $\beta$ and any $c \in \res^\times$, we define 
\begin{equation}\label{E:nh}
n_\beta(c) = u_\beta(c)u_{-\beta}(-c^{-1})u_\beta(c) \quad\quad \text{and} \quad\quad
 h_{\beta}(c) = n_\beta(c) n_\beta(1)^{-1}.
 \end{equation} 
 For further properties of the elements $n_\beta(c)$ and $h_\beta(c)$, we refer the reader to~\cite[Section~3]{PRS} and the references therein.

In order to prove Proposition~\ref{prop:doublecoset} below, we will make use of three important relations in $G(F)$.  We will see in the proof of Proposition~\ref{prop:doublecoset} that these relations allow us to rewrite a point in $IxI$ as a point in $(I_P)^y z I$ for certain $z \in W$.  The fact these these relations exist follows from \cite{BilligDyer}, for instance.

\begin{enumerate}
\item For all affine roots $\beta$, all $b \in I$, and all $c \in \res$, there exists a unique $b' \in I$ and a unique $\tilde{c} \in \res$ such that
\begin{equation}\label{eq:coset}
b u_\beta(c)s_\beta = u_\beta(\tilde{c})s_\beta b'.
\end{equation}

\item For all affine roots $\beta, \beta'$ and all $c \in \res$,
\begin{equation}\label{eq:fconj}
s_\beta u_{\beta'}(c)s_\beta= u_{s_{\beta}(\beta')}(\pm c), 
\end{equation} 
where the sign for $\pm c$ is uniquely determined by the pair $\beta, \beta'$.

\item 
For all affine roots $\beta$ and all $c \in \res^\times$, 
\begin{equation} \label{eq:MFL}
u_\beta(c)s_\beta = u_{-\beta}(c^{-1})u_\beta(-c)h_{\beta}(c). \quad\quad \mbox{(Main Folding Law)}
\end{equation}
\end{enumerate}

We now use the above relations to connect certain double cosets in $G(F)$ to labeled folded alcove walks which are positively folded with respect to a $P$-chimney. The structure of the induction and resulting cases in the proof of Proposition \ref{prop:doublecoset} below are modeled upon similar results in \cite{PRS, MR}.

\begin{prop} \label{prop:doublecoset} Let $P = MN$ be a standard spherical parabolic subgroup of $G(F)$.
	Then:
	\begin{enumerate}
		\item For all $x \in \aW$, every point of $IxI$ is a point of $I_P z I$ for some $z \in \aW$.
		\item For all $x, z \in \aW$, the points of $IxI \cap I_P z I$ are in bijection with the set of labeled folded alcove walks of type $\vec{x}$ from $\fa$ to $\z$ which are positively folded with respect to the $P$-chimney. 
	\end{enumerate}
\end{prop}

\begin{proof} The proof is by induction on the length of $x$. We begin by considering the case $\ell(x)=1$. In this case, $x=s_j$ for some $j \in \{0,1,\dots, n\}$.  Recall that for $j \in \{1,\dots,n\}$, the spherical root $\alpha_j \in \Delta$ is a positive simple root with respect to $B^-$, while $\alpha_0$ is the affine root such that $U_{\alpha_0} = U_{-\widetilde\alpha,-1}$.  Note also that the subgroup $U_{\alpha_j}$ is contained in $I$ for all $j \in \{0,1,\dots,n\}$.  If $j \in \{1,\dots,n\}$, this follows from $\alpha_j$ being positive with respect to $B^-$ and $I$ being the preimage of $B^-(\res)$.  For $j = 0$ recall that by definition $U_{\alpha_0}$ fixes the half-apartment bounded by $H_{\widetilde\alpha,1}$ which contains $\fa$, hence is contained in $I$.  Similarly, for all $j \in \{0, 1, \dots, n\}$, we have that $\overline U_{\alpha_j}$ is a subgroup of $I$.
	
	Now by Theorem 6.15 in \cite{Ronan}, for example, we can write the points of $IxI = I s_j I$ as 
	\begin{equation}\label{E:IsI}
	Is_jI=\left\{u_{\alpha_j}(c)s_j I \mid c \in \res\right\}.
	\end{equation}
Observe also that there is a unique minimal gallery $\gamma_{s_j} = (p_0, c_0, p_1, c_1,p_2)$ of type $\vec{s}_j$ with first face $p_0 = c_0 = \fa$.  This gallery $\gamma_{s_j}$ has final face $c_1 = p_2 = s_j \fa$, and the supporting hyperplane of its panel $p_1$ is $H_{\alpha_j,0}$ if $j \in \{1,\dots,n\}$ and $H_{\widetilde \alpha,1}$ if $j  = 0$.
	
	We first suppose that $j \in \{1,\dots,n\}$ and that $U_{\alpha_j}$ is contained in $M$.  Equivalently, the positive simple root $\alpha_j$ lies in $\Phi^+_M$.  In this case, the minimal gallery $\gamma_{s_j}$ crosses $p_1 \subset H_{\alpha_j,0}$ from the negative to the positive side.  Thus there is no positive fold possible at $p_1$, and so all positively folded galleries of type~$\vec{s}_j$ with first face $\fa$ are in fact equal to $\gamma_{s_j}$.  Hence all such galleries have final alcove~$s_j\fa$.  It follows that for $z \neq s_j$, there are no labeled folded alcove walks of type $\vec{s}_j$ from $\fa$ to $\z$ which are positively folded with respect to the $P$-chimney.  If $z = s_j$, then by labeling $p_1$ by any element of $\res$, we obtain a labeled folded alcove walk of type $\vec{x} = \vec{s}_j$ from $\fa$ to $\z = s_j\fa$, and each such labeling of $p_1$ corresponds to a distinct labeled folded alcove walk.  
	
	We then observe that $u_{\alpha_j}(c) \in I \cap M$ for any $c \in \res$, so $u_{\alpha_j}(c) \in I_P$ and  thus $u_{\alpha_j}(c)s_j I \in I_P s_j I$.  Combining this with~\eqref{E:IsI} we obtain that every point of $Is_jI$ is a point of $I_P s_j I$, which establishes (1) in this case.  Using Corollary~\ref{cor:disjoint}, it follows that for all $z \neq s_j$, the intersection $Is_j I \cap I_P z I$ is empty.  For $z = s_j$, the map taking the element $c \in \res$ which labels the panel $p_1$ of $\gamma_{s_j}$  to $u_{\alpha_j}(c)s_j I$ induces a bijection between the set of labeled folded alcove walks of type $\vec{x} = \vec{s}_j$ from $\fa$ to $\z = s_j\fa$ which are positively folded with respect to the $P$-chimney, and the points of $I s_j I \cap I_P s_j I$.  This completes the proof of (2) in this case.
	
	Next we suppose that $j \in \{1,\dots,n\}$ and that $U_{\alpha_j}$ is not contained in $M$ (hence $U_{-\alpha_j}$ will be contained in $N$).  In other words, $\alpha_j \in \Phi^+ \setminus \Phi^+_M$.  In this case, the gallery $\gamma_{s_j}$ crosses $p_1 \subset H_{\alpha_j,0}$ from the positive to the negative side.  We thus have two possible positively folded galleries of type $\vec{x}= \vec{s}_j$ with first face $\fa$: one equal to $\gamma_{s_j}$, and the other obtained by folding $\gamma_{s_j}$ in the panel $p_1$, so as to have final face~$\fa$.  Thus for $z \not \in \{ s_j, \id \}$, there are no labeled folded alcove walks of type $\vec{s}_j$ from $\fa$ to $\z$ which are positively folded with respect to the $P$-chimney.  If $z = s_j$ then the panel $p_1$ can only be labeled by $0$, while if $z = \id$ then the panel $p_1$ can be labeled by any element of $\res^\times$.
	
	Now consider the point $u_{\alpha_j}(c) s_j I$ of $Is_jI$.  If $c = 0$, then $u_{\alpha_j}(c)$ is trivial so obviously $u_{\alpha_j}(c)s_j I \in I_P s_j I$.  However, if $c \neq 0$, then $u_{\alpha_j}(c)\in I$ but $u_{\alpha_j}(c)\notin I_P$.  By the Main Folding Law (Equation \eqref{eq:MFL}), for all $c \neq 0$
	$$u_{\alpha_j}(c)s_j = u_{-\alpha_j}(c^{-1})u_{\alpha_j}(-c)h_{\alpha_j}(c).$$
	Since $\alpha_j \notin R(I_P)$, we have $-\alpha_j \in R(I_P)$, and so $u_{-\alpha_j}(c^{-1}) \in I_P$. Then since $u_{\alpha_j}(-c)$ and $h_{\alpha_j}(c)$ are both elements of $I$, we obtain that $u_{\alpha_j}(c)s_j I =u_{-\alpha_j}(c^{-1})u_{\alpha_j}(-c)h_{\alpha_j}(c)I \in I_P \id I$ for all $c \neq 0$.  We have now proved (1) in this case.  We also see, using Corollary~\ref{cor:disjoint}, that for all $z \not \in \{ s_j, \id \}$, the intersection $I s_j I \cap I_P z I$ is empty.
	
	For (2), the previous paragraph shows that the unique point of $Is_jI \cap I_P s_j I$ is $u_{\alpha_j}(0)s_j I = s_jI$.  Letting $z = s_j$, we map the (unique) labeled folded alcove walk of type $\vec{x} = \vec{s}_j$ from $\fa$ to $\z = s_j \fa$ with panel $p_1$ labeled by $0$ to the (unique) point $u_{\alpha_j}(0)s_j I$ of $I s_j I \cap I_P s_j I$.  Now suppose $c \neq 0$.  Then for $z = \id$, the map which takes the label $c \in \res^\times$ of the panel $p_1$ in the gallery from $\fa$ to $\fa$ obtained by folding $\gamma_{s_j}$ in $H_{\alpha_j,0}$ to the coset $u_{\alpha_j}(c)s_j I$ induces a bijection between the set of labeled folded alcove walks of type $\vec{x} = \vec{s}_j$ from $\fa$ to $\z = \fa$ which are positively folded with respect to the $P$-chimney, and the points of $I s_j I \cap I_P z I = Is_j I \cap I_P \id I$.  This completes the proof of (2) in this case.
	
	Finally, we consider $j = 0$; that is, $u_{\alpha_0}(c)s_0 I \in I s_0 I$.  In this case, $\gamma_{s_0}$ crosses $p_1 \subset H_{\widetilde \alpha,1}$ from the negative to the positive side, and $u_{\alpha_0}(c) \in I \cap I_P$ for all $c\in \res$.  The remainder of the proofs of (1) and (2) in this case are then similar to the first case above, in which $j \in \{1,\dots,n\}$ and $\alpha_j \in \Phi_M^+$.  We have now established the base case of our induction.
	
	For the inductive step, let $(s_{j_1}, s_{j_2}, \dots, s_{j_l})$ be a reduced word for $x$, and let $s_j \in \aS$ be such that $\ell(xs_j)>\ell(x)$.  Using Theorem 6.15 and then Lemma 7.4 of \cite{Ronan}, we consider the element
	$u_{j_1}(c_1)s_{j_1} \dots u_{j_l}(c_l)s_{j_l} u_{\alpha_j}(c)s_j$ of $Ixs_j I,$ where $c_i, c \in \res$, and we have put $u_{j_i}(c)$ for $u_{\alpha_{j_i}}(c)$.  Noting that $u_{j_1}(c_1)s_{j_1} \dots u_{j_l}(c_l)s_{j_l}$ is a point of $IxI$, we have by  the inductive hypothesis for~(1) that there exists  $\tilde{u} \in I_P$, $z \in \aW$, and $b \in I$ such that  
	\begin{equation}\label{eqn:induction} \left[ u_{j_1}(c_1)s_{j_1} \dots u_{j_l}(c_l)s_{j_l}\right] u_{\alpha_j}(c)s_j = \left[\tilde{u} z b \right] u_{\alpha_j}(c)s_j. \end{equation}
	Now by Equation \eqref{eq:coset} there exist unique $b' \in I$ and $\tilde{c} \in \res$ such that $b u_{\alpha_j}(c)s_j=u_{\alpha_j}(\tilde{c})s_j b'$, and so
	$$ \tilde{u}zb u_{\alpha_j}(c)s_j = \tilde{u}z u_{\alpha_j}(\tilde{c})s_j b'.$$ 
	
	We now finish the proof of (1) by considering two cases.
	
	\medskip
	\noindent \emph{Case 1:} The affine root $z(\alpha_j)$ is in $R(I_P)$. 
	
	In this case, we observe that we have the equality 
	$$\tilde{u}z u_{\alpha_j}(\tilde{c}) s_j b'= \tilde{u}z u_{\alpha_j}(\tilde{c})z^{-1}z s_j b',$$
	and then by Equation \eqref{eq:fconj}, we have
	$$\tilde{u}z u_{\alpha_j}(\tilde{c})z^{-1}z s_j b'= \tilde{u}u_{z(\alpha_j)}(\pm \tilde{c})z s_j b'.$$
	Since $z(\alpha_j) \in R(I_P)$, 
	\begin{equation}\label{eq:Case1} \tilde{u}u_{z(\alpha_j)}(\pm \tilde{c})z s_j b' \in I_P zs_j I.
	\end{equation}
	Combining this with Equation~\eqref{eqn:induction} we obtain that for any $c_i, c \in \res$, 
	\begin{equation}u_{j_1}(c_1)s_{j_1} \dots u_{j_l}(c_l)s_{j_l} u_{\alpha_j}(c)s_j \in I_P zs_j I.
	\end{equation}
	This proves (1) in this case.
	
	\medskip
	\noindent \emph{Case 2:} The affine root $z(\alpha_j)$ is not in $R(I_P)$.  
	
	We will consider two subcases.
	
	\medskip\noindent \emph{Case 2(a):} $z(\alpha_j) \notin R(I_P)$ and $\tilde{c} \neq 0$. 
	
	By the Main Folding Law (Equation \eqref{eq:MFL}), we have 
	\begin{eqnarray*}
		\tilde{u}z u_{\alpha_j}(\tilde{c}) s_j b' &=& \tilde{u}z u_{-\alpha_j}(\tilde{c}^{-1})u_{\alpha_j}(-\tilde{c})h_{\alpha_j}(\tilde{c}) b' \\
		&=&\tilde{u}z u_{-\alpha_j}(\tilde{c}^{-1})b'' 
	\end{eqnarray*}
	where $b''=u_{\alpha_j}(-\tilde{c})h_{\alpha_j}(\tilde{c}) b' \in I$. By Equation \eqref{eq:fconj}, we have
	\begin{eqnarray*}
		\tilde{u}z u_{-\alpha_j}(\tilde{c}^{-1})b'' &=&\tilde{u}z u_{-\alpha_j}(\tilde{c}^{-1})z^{-1}z b'' \\
		&= &\tilde{u}u_{z(-\alpha_j)}(\pm \tilde{c}^{-1})z b''.
	\end{eqnarray*} 
	Since $u_{z(-\alpha_j)}(\pm \tilde{c}^{-1}) \in I_P$, 
	\begin{equation}\label{eq:Case2a} \tilde{u}u_{z(-\alpha_j)}(\pm \tilde{c}^{-1})z b'' \in I_P z I. \end{equation} 
	Using Equation~\eqref{eqn:induction} we get that
	\begin{equation}u_{j_1}(c_1)s_{j_1} \dots u_{j_l}(c_l)s_{j_l} u_{\alpha_j}(c)s_j \in I_P z I. \end{equation} 
	This proves (1) in this case.
	
	\medskip\noindent  \emph{Case 2(b):} $z(\alpha_j) \notin R(I_P)$ and $\tilde{c} = 0$.  
	
	We observe that $u_{\alpha_j}(0)=u_{-\alpha_j}(0)$, and so we have 
	\begin{eqnarray*}
		\tilde{u}z u_{\alpha_j}(0)z^{-1}z s_j b' &=& \tilde{u}z u_{-\alpha_j}(0)z^{-1}z s_j b' \\
		&=& \tilde{u}u_{z(-\alpha_j)}(0) z s_j b'.
	\end{eqnarray*}
	Since $u_{z(-\alpha_j)}\in {I_P},$ 
	\begin{equation}\label{eq:Case2b} \tilde{u}u_{z(-\alpha_j)}(0) z s_j b' \in I_P zs_j I. \end{equation}
	Using Equation~\eqref{eqn:induction} we get that 
	\begin{equation} u_{j_1}(c_1)s_{j_1} \dots u_{j_l}(c_l)s_{j_l} u_{\alpha_j}(c)s_j \in I_P zs_j I. \end{equation}
	This proves (1) in this case, and so completes the proof of the inductive step for (1).
	
	\medskip
	
	In order to complete the proof of (2), we first observe that by combining Equations~\eqref{eqn:induction}, \eqref{eq:Case1},~\eqref{eq:Case2a}, and~\eqref{eq:Case2b}, we have that whenever $\ell(x s_j) > \ell(x)$, every point of $Ixs_jI$ is contained in either $I_P z I$ or $I_P zs_j I$ for some $z \in \aW$ such that $IxI \cap I_P zI$ is nonempty.  Hence by Corollary~\ref{cor:disjoint}, the only double coset intersections $Ixs_j I \cap I_P z' I$ which are nonempty are those where $z' \in \{ z, zs_j\}$ and $IxI \cap I_P z I$ is nonempty. 
	
	Now suppose there exists a labeled folded alcove walk of type $\overrightarrow{x s_j}$ from $\fa$ to some alcove $z'\fa$ which is positively folded with respect to the $P$-chimney.  As $\ell(xs_j) > \ell(x)$, we may assume that the last panel of this alcove walk is of type $s_j$.  Since there may be either a fold or a crossing at this last panel, by truncation we obtain a labeled folded alcove walk of type $\vec{x}$ from $\fa$ to either $z'\fa$ or $z' s_j \fa$, respectively, which is positively folded with respect to the $P$-chimney.  If $IxI \cap I_P z' I$ (respectively, $IxI \cap I_P z' s_j I$) is empty, this contradicts the inductive assumption for (2).  It follows that all labeled folded alcove walks of type $\overrightarrow{x s_j}$ with first alcove $\fa$ which are positively folded with respect to the $P$-chimney have final alcove either $z\fa$ or $z s_j\fa$, where $z \in \aW$ is such that $IxI \cap I_P z I$ is nonempty.  
	
	Therefore to complete the proof of (2), we just need to consider the points of $Ixs_jI \cap I_P z' I$ where $z' \in \{ z, zs_j\}$ and $IxI \cap I_P z I$ is nonempty.  By the inductive hypothesis for (2), each point of $IxI \cap I_P z I$ corresponds to a labeled folded alcove walk of type $\vec{x}$ from $\fa$ to $\z$ which is positively folded with respect to the $P$-chimney.  Since $\ell(x s_j) > \ell(x)$, a gallery of type $\overrightarrow{xs_j}$ starting at $\fa$ can be obtained from a gallery of type $\vec{x}$ from $\fa$ to $\z$ by adding either a crossing or a fold at the panel of $\z$ supported by $H_{z(\alpha_j)}$. 
	
	We will use the same cases as above to complete the proof. 
	
	\medskip
	\noindent \emph{Case 1:} The affine root $z(\alpha_j)$ is in $R(I_P)$. 
	
	A gallery of type $\overrightarrow{xs_j}$ from $\fa$ to $zs_j \fa$ can be obtained from a gallery of type $\vec{x}$ from $\fa$ to $\z$ by adding a crossing from $z\fa$ to $zs_j\fa$.  In this case, crossing $H_{z(\alpha_j)}$ from $z\fa$ to $zs_j\fa$ goes from the negative to the positive side.  Thus the panel of $H_{z(\alpha_j)}$ at which this crossing occurs can be labeled by any element of $\res$.  
	Now use the element $\pm \tilde{c} \in \res$ from Equation~\eqref{eq:Case1} as this label.  By the inductive hypothesis for (2), this induces a bijection between the set of labeled folded alcove walks of type $\overrightarrow{xs_j}$ from $\fa$ to $zs_j\fa$ which are positively folded with respect to the $P$-chimney, and the points of $I x s_j I \cap I_P zs_j I$.
	
	\medskip
	\noindent \emph{Case 2:} The affine root $z(\alpha_j)$ is not in $R(I_P)$.  
	
	In this case, a gallery crossing the hyperplane $H_{z(\alpha_j)}$ from $z\fa$ to $zs_j\fa$ goes from the positive to the negative side.  Hence a gallery $\gamma$ of type $\overrightarrow{xs_j}$ with initial subgallery of type $\vec{x}$ going from $\fa$ to $z\fa$ can either be positively folded in $H_{z(\alpha_j)}$ or cross $H_{z(\alpha_j)}$. Folding in $H_{z(\alpha_j)}$ implies that $\gamma$ ends at $z\fa$ and the panel of $H_{z(\alpha_j)}$ at which this fold occurs can be labeled by any element of $\res^\times$. If $\gamma$ has a crossing at $H_{z(\alpha_j)}$ it will end at $zs_j\fa$, with the panel of $H_{z(\alpha_j)}$ at which this crossing occurs labeled by $0$.  We will see that the subcases below correspond to these two possibilities, respectively.
	
	\medskip\noindent \emph{Case 2(a):} $z(\alpha_j) \notin R(I_P)$ and $\tilde{c} \neq 0$. 
	
	Here, we can use the element $\pm \tilde{c}^{-1} \in \res^\times$ from Equation~\eqref{eq:Case2a} to label the panel of $H_{z(\alpha_j)}$ at which a gallery ending at $z\fa$ (as described in the first possibility above) has a positive fold.  By the inductive hypothesis for (2), this induces a bijection between the set of labeled folded alcove walks of type $\overrightarrow{xs_j}$ from $\fa$ to $z\fa$ which are positively folded with respect to the $P$-chimney, and the points of $I x s_j I \cap I_P z I$. 
	
	\medskip\noindent  \emph{Case 2(b):} $z(\alpha_j) \notin R(I_P)$ and $\tilde{c} = 0$.  
	
	Now we use $0$ to label the panel of $H_{z(\alpha_j)}$ between $z\fa$ and $zs_j\fa$, as in the second possibility described above. By Equation~\eqref{eq:Case2b} and the inductive hypothesis for (2), this induces a bijection between the set of labeled folded alcove walks of type $\overrightarrow{xs_j}$ from $\fa$ to $zs_j\fa$ which are positively folded with respect to the $P$-chimney, and the points of $Ixs_jI \cap Izs_jI$.
	
	\medskip
	
	This completes the proof of Proposition~\ref{prop:doublecoset}.
\end{proof}

The next result extends Proposition~\ref{prop:doublecoset}(2) to $(P,y)$-chimneys.

\begin{thm} \label{thm:doublecosety}  Let $x, y, z \in \aW$ and let $P$ be a standard spherical parabolic subgroup of $G(F)$.  Then there is a bijection between the points of the intersection 
	\[  I x I  \cap (I_P)^y z I \] 
	and the set of labeled folded alcove walks of type $\vec{x}$ from $\fa$ to $\z$ which are positively folded with respect to the $(P,y)$-chimney.
\end{thm}

\begin{proof} Suppose $u_\beta(c) \in (I_P)^y$ for some affine root $\beta$ and $c \in \res$. This means that $y^{-1}u_\beta(c)y \in I_P$.  By Equation \eqref{eq:fconj}, $y^{-1}u_\beta(c)y= u_{y^{-1}(\beta)}(\pm c)$. Therefore, the result follows from Proposition~\ref{prop:doublecoset}(2) after transforming the orientation via the left action of $y$ as described in \cite[Definition  3.4]{GraeberSchwer}. 
\end{proof}

\begin{remark} \label{rem:braidDoubleCoset} We note that Theorem \ref{thm:doublecosety} does not depend on a choice of minimal word for~$x$. Given two different reduced words for $x$, we can form minimal galleries of distinct types $\vec{x}$ and $\vec{x}'$ corresponding to these words. For any $z\in W$ such that there exists a $(P,y)$-positively folded gallery of type $\vec{x}$ from $\fa$ to $\z$, we know by Theorem \ref{thm:doublecosety} that the intersection $I x I  \cap (I_P)^y z I$ is nonempty. The theorem further implies that if $I x I  \cap (I_P)^y z I$ is nonempty, there exists a $(P,y)$-positively folded gallery of type $\vec{x}'$ from $\fa$ to $\z$, since $\vec{x}'$ corresponds to another reduced expression for $x$. Thus, we have yet another proof of the braid-invariance of $(P,y)$-chimney orientations (compare Proposition~\ref{prop:braidinvariant} and Remark~\ref{rem:braidRecursion}).
\end{remark}

We now use Theorem~\ref{thm:doublecosety} to prove Theorem \ref{thm:DoubleCosetIntersectionParahoric}, the main result of this section. 
Here, for any face $\sigma$ of $\fa$ which contains the origin $v_0$, we denote by $K(\sigma)$ the parahoric subgroup of $G(F)$ which is the stabilizer in $G(F)$ of $\sigma$. (To avoid confusion with the spherical parabolic subgroup $P$, we are using the letter $K$ for parahorics.)  Thus in particular, $K = K(v_0)$ and $I = K(\fa)$.  We write $W_\sigma$ for the stabilizer of $\sigma$ in $\sW$.

\begin{definition} Let $\sigma$ and $\tau$ be faces of $\fa$ which contain the origin.  Then $x\in W$ is:
	\begin{enumerate}
		\item \textit{left-$W_\sigma$-reduced} if $\ell(wx)\geq \ell(x)$ for all $w \in W_\sigma$;
		\item \textit{right-$W_\tau$-reduced} if $\ell(xv)\geq \ell(x)$ for all $v \in W_\tau$; and
		\item \textit{$(W_\sigma,W_\tau)$-reduced} if $x$ is left-$W_\sigma$-reduced and right-$W_\tau$-reduced, i.e. $\ell(wx)\geq \ell(x)$ for all $w \in W_\sigma$ and $\ell(xv)\geq \ell(x)$ for all $v \in W_\tau.$
	\end{enumerate} 
\end{definition}

\begin{thm}\label{thm:DoubleCosetIntersectionParahoric} Let $P$ be a standard spherical parabolic subgroup of $G(F)$. For any faces $\sigma$ and $\tau$ of $\fa$ which contain the origin and any $x, y, z \in \aW$ such that $x$ is $(W_\sigma,W_\tau)$-reduced,
	there is a bijection between the points of the intersection 
	\[  K(\sigma) x K(\tau)  \cap (I_P)^y z K(\tau) \] 
	and the union over $w\in W_\sigma$ of the set of labeled folded alcove walks of type $\overrightarrow{wx}$ from $\fa$ to $\z$ which are positively folded with respect to the $(P,y)$-chimney. 
\end{thm}

\noindent Note that the hypothesis that $x$ is $(W_\sigma,W_\tau)$-reduced does not lead to a loss of generality, since we can always find a representative of the double coset $K(\sigma) x K(\tau)$ that is of this form. This hypothesis is included only so that we can more easily describe the type of labeled walks required.  The proof of this theorem uses techniques inspired by the proof of the equivalence of Theorems G and C in \cite{ParkinsonRam}. 

\begin{proof}  We first consider the case in which $\tau = \fa$, so that $K(\tau)=I$. Recall that $K(\sigma)=\sqcup_{w \in W_\sigma} IwI$.   Since $x$ is left-$W_\sigma$-reduced, for any $w\in W_\sigma$ we have $(IwI)xI=IwxI$. This implies that \[ K(\sigma) x I=\bigsqcup_{w \in W_\sigma} (IwI) x I=\bigsqcup_{w \in W_\sigma} IwxI.\] Thus, by Theorem~\ref{thm:doublecosety}, the points of $K(\sigma) x I  \cap (I_P)^y z I$ are in bijection with the (disjoint) union over $w \in W_\sigma$ of the set of labeled folded alcove walks of type $\overrightarrow{wx}$ from $\fa$ to $\z$ which are positively folded with respect to the $(P,y)$-chimney. This completes the proof in the case that $\tau = \fa$.
	
	Now suppose $\sigma=\fa$ and let $\tau$ be any face of $\fa$ which contains the origin.  We may assume without loss of generality that $\tau \neq \fa$.  Fix a minimal gallery $\gamma$ in the standard apartment $\App$ from $\fa$ to $x\tau$.  The points of $I x K(\tau)$ are then in bijection with the set of minimal galleries in the building $X$ which are of the same type as $\gamma$.  
	
	Denote by $\gamma'$ the gallery in $\App$ which is obtained from $\gamma$ by removing its final simplex $x\tau$.  We claim that there is a bijection between the set of minimal galleries in $X$ of the same type as $\gamma$ and the set of minimal galleries in $X$ of the same type as $\gamma'$.  To see this, note that the set of alcoves in $\App$ which contain the face $\tau$ of $\fa$ is given by $\{ v\fa \mid v\in W_\tau \}$.  Hence the set of alcoves in $\App$ which contain the face $x\tau$ of $\x=x\fa$ is given by $\{ xv\fa \mid v\in W_\tau \}$.  Since $x$ is right-$W_\tau$-reduced and $\gamma$ is a minimal gallery from $\fa$ to $x\tau$, it follows that $\gamma$ must have final alcove $\x$, and so the truncated gallery $\gamma'$ is a minimal gallery from $\fa$ to $\x$.  A similar  truncation process establishes the claimed bijection between sets of minimal galleries in $X$.
	
	Since the points of $IxI$ are in bijection with the set of minimal galleries in $X$ which are of the same type as $\gamma'$, we now have a bijection between the points of $I x K(\tau)$ and those of $IxI$.  By Theorem~\ref{thm:doublecosety}, the points of the intersection $I x I \cap (I_P)^y z I$ are in bijection with the set of labeled folded alcove walks of type $\vec{x}$ from $\fa$ to $\z$ which are positively folded with respect to the $(P,y)$-chimney.  Projecting points of $(I_P)^y z I$ back to the flag variety $G/K(\tau)$ gives us points in $(I_P)^y z K(\tau)$. Thus, points in $I x K(\tau)  \cap (I_P)^y z K(\tau)$ are also in bijection with the set of labeled folded alcove walks of type $\vec{x}$ from $\fa$ to $\z$ which are positively folded with respect to the $(P,y)$-chimney.  This completes the proof in the case that $\sigma = \fa$.
	
	Finally, let $\sigma$ and $\tau$ be any faces of $\fa$ which contain the origin. As before, we observe that since $x$ is left-$W_\sigma$-reduced,
	\[ K(\sigma) x K(\tau)=\bigsqcup_{w \in W_\sigma} (IwI) x K(\tau)=\bigsqcup_{w \in W_\sigma} IwxK(\tau).\]
	Also, since $x$ is right-$W_\tau$-reduced, the points of $K(\sigma) x K(\tau) \cap (I_P)^y z K(\tau)$ are in bijection with the points of $K(\sigma) x I \cap (I_P)^y z I$. Thus in general, the points of $K(\sigma) x K(\tau) \cap (I_P)^y z K(\tau)$ are in bijection with the union over $w \in W_\sigma$ of the set of labeled folded alcove walks of type $\overrightarrow{wx}$ from $\fa$ to $\z$ which are positively folded with respect to the $(P,y)$-chimney.   
\end{proof}

\begin{remark}\label{rem:APVM}
	In \cite{APvM}, Abramenko, Parkinson, and Van Maldeghem study and count the number of points in intersections of spheres (for the Weyl-distance) in certain buildings. While we assume for all of Section~\ref{sec:galleriesCosets} that the building in question is a Bruhat--Tits building, and is thus of affine type, Abramenko, Parkinson, and Van Maldeghem study buildings of arbitrary type, but assume in addition that they are locally finite and regular. There is, however, some similarity to our results. Theorem 4.7 in \cite{APvM} for example may be translated into a statement about the cardinality of the double coset intersections appearing in the statement of Theorem~\ref{thm:DoubleCosetIntersectionParahoric}.  
\end{remark}

The next result determines nonemptiness of the double coset intersections appearing in Theorem~\ref{thm:DoubleCosetIntersectionParahoric}.

\begin{thm}\label{thm:IntersectionParahoricShadow} Let $P$ be a standard spherical parabolic subgroup of $G(F)$. For any faces $\sigma$ and $\tau$ of $\fa$ which contain the origin and any $x, y, z \in \aW$ such that $x$ is $(W_\sigma,W_\tau)$-reduced, the intersection 
	\[  K(\sigma) x K(\tau)  \cap (I_P)^y z K(\tau) \] 
	is nonempty if and only if there exists a gallery of type $(\type(\sigma),\vec{x},\type(\tau))$ from $\sigma$ to $z\tau$ which is positively folded with respect to the $(P,y)$-chimney. Equivalently, \[  K(\sigma) x K(\tau)  \cap (I_P)^y z K(\tau) \neq \emptyset \] if and only if the simplex $z\tau$ lies in the shadow $\Sh_{P,y}(x\tau,\sigma)$. 
\end{thm}

\noindent As with Theorem~\ref{thm:DoubleCosetIntersectionParahoric}, the hypothesis on $x$ here does not lead to any loss of generality.

\begin{proof} Suppose first that $K(\sigma) x K(\tau)  \cap (I_P)^y z K(\tau) \neq \emptyset$. Then by Theorem \ref{thm:DoubleCosetIntersectionParahoric}, for some $w \in W_\sigma$ there is a labeled folded alcove walk $\gamma$ of type $\overrightarrow{wx}$ from $\fa$ to $\z$ which is positively folded with respect to the $(P,y)$-chimney.  We now use $\gamma$ to construct a gallery $\gamma''$ of type $(\type(\sigma),\vec{x},\type(\tau))$ from $\sigma$ to $z\tau$ which is positively folded with respect to the $(P,y)$-chimney.  Since $x$ is $W_\sigma$-reduced, we may choose a reduced word $(s_{i_1}, \dots, s_{i_r})$ for $wx$ such that $s_{i_1}\dots s_{i_t} = w$ and $s_{i_{t+1}}\dots s_{i_r}=x$, and assume that $\gamma$ is of type this word.  Letting $\gamma = (p_0, c_0, p_1, c_1, \dots, p_r, c_r,c_{r+1})$, so that $p_0 = c_0 = \fa$ and $c_r = p_{r+1} = \z$, we may thus assume that the panel $p_{j}$ of $\gamma$ has type $s_{i_j}$ for $1 \leq j \leq r$.  This means that the first $t$ panels $p_1,\dots,p_t$ of $\gamma$ must each contain $\sigma$, and so the alcove $c_t$ of $\gamma$ lies in the star of $\sigma$ in $\App$.  Therefore $c_t=w'\fa$ for some $w' \in W_\sigma$.  
	
	Now define $\gamma'$ to be the subgallery of $\gamma$ which starts at $c_t$ and continues in the same way as $\gamma$ to $\z$.  Then $\gamma'$  is a gallery of type $\vec{x}$ from $\w'$ to $\z$.  Finally let $\gamma''$ be the gallery obtained from $\gamma'$ by changing its first simplex to $\sigma$ and its final simplex to $z\tau$.  Then by our choice of word for $wx$ and $x$ being right-$W_\tau$-reduced, the gallery $\gamma''$ from $\sigma$ to $z\tau$ is of type $(\type(\sigma),\vec{x},\type(\tau))$.  All folds in $\gamma''$ are folds of the original gallery $\gamma$, and so $\gamma''$ is a gallery of the desired type which is positively folded with respect to the $(P,y)$-chimney.
	
	Conversely, suppose we have a gallery $\gamma''$ of type $(\type(\sigma),\vec{x},\type(\tau))$ from $\sigma$ to $z\tau$ which is positively folded with respect to the $(P,y)$-chimney.  Let $\w$ be the first alcove of $\gamma''$ and note that $\w$ contains $\sigma$.  Since $x$ is left-$W_\sigma$-reduced, the second alcove of $\gamma''$ cannot contain $\sigma$.  Similarly, the final alcove of $\gamma''$ contains $z\tau$, and since $x$ is right-$W_\tau$-reduced, the final alcove of $\gamma''$ must be $\z$.  Now define $\gamma$ to be the gallery obtained from $\gamma''$ by replacing its initial simplex $\sigma$ by a minimal gallery from $\fa$ to $\w$, and replacing its final simplex $z\tau$ by $\z$.  Then $\gamma$ is a gallery of type $\overrightarrow{wx}$ from $\fa$ to $\z$, for a suitable choice of reduced word for $wx$.  By construction the initial subgallery of $\gamma$ from $\fa$ to $\w$ is minimal hence has no folds, while any folds in the portion of $\gamma$ from $\w$ to $\z$ are also folds of $\gamma''$, thus $\gamma$ is  positively folded with respect to the $(P,y)$-chimney.
	
	We can then obtain a labeled folded alcove walk of type $\overrightarrow{wx}$ from $\fa$ to $\z$ which is positively folded with respect to the $(P,y)$-chimney by choosing a labeling each of the panels of $\gamma$ according to the scheme given by Definition~\ref{def:labeledwalks}.  Any such labeled folded alcove walk corresponds to a point of the intersection $K(\sigma) x K(\tau)  \cap (I_P)^y z K(\tau)$ by Theorem \ref{thm:DoubleCosetIntersectionParahoric}. Therefore, this intersection is nonempty.
\end{proof}

\begin{remark} In light of the above results, it is tempting to try to formulate a bijection directly between the points of the intersection $K(\sigma) x K(\tau)  \cap (I_P)^y z K(\tau)$ and suitably labeled galleries of type $(\type(\sigma),\vec{x},\type(\tau))$ from $\sigma$ to $z\tau$ which are positively folded with respect to the $(P,y)$-chimney. However, the following example shows that this is not always a useful compression.
	
	In type $\tilde{A}_2$, let $\tau=\fa$, so that $K(\tau)=I$, and let $\sigma=v_0$, so that $K(\sigma)=K$. Let $x=s_0 \in \aS$ be the reflection in the (affine) hyperplane $H_{\widetilde\alpha,1}$ which bounds $\fa$.  Then \[K(\sigma)x K(\tau) = Ks_0 I=\bigsqcup_{w\in W_0} Iws_0I.\] In particular, the cells $Is_0I$, $Is_1s_0I$, and $Is_2s_0I$ are mutually disjoint.  Now let $P=B$ and $y = \id$, so that the $(P,y)$-chimney is the $B$-chimney, and let $z = s_0$.  We claim that each of these three cells has nonempty intersection with $(I_P)^y z I = I_P s_0 I$.  
	
	By Theorem~\ref{thm:doublecosety}, the intersections of $Is_0I$, $Is_1s_0I$, and $Is_2s_0I$ with $I_P s_0 I$ are nonempty if and only if there exist labeled folded alcove walks from $\fa$ to $s_0\fa$ of types $\overrightarrow{s_0}$, $\overrightarrow{s_1s_0}$, and $\overrightarrow{s_2s_0}$, respectively.  Consider the galleries $(\fa, p_0, s_0\fa)$, $(\fa, p_1,  \fa, p_0, s_0\fa)$, and $(\fa, p_2,\fa,p_0,s_0\fa)$ from $\fa$ to $s_0\fa$, where for simplicity we have omitted the first face $\fa$ and last face $s_0\fa$ of each gallery, and for $i \in \{0,1,2\}$ we abuse notation and write $p_i$ for the panel of $\fa$ of type~$s_i$.  These galleries have types $\overrightarrow{s_0}$, $\overrightarrow{s_1s_0}$, and $\overrightarrow{s_2s_0}$, respectively, and are positively folded with respect to the $B$-chimney.  Thus any labeling of their panels corresponds to a point in the intersection of $Is_0I$, $Is_1s_0I$, and $Is_2s_0I$ with $I_P s_0 I$, respectively.  This proves the claim.
	
	Hence if we wished to instead enumerate points in $K s_0 I \cap I_P s_0 I$ by labeling just a gallery of the form $(v_0, \fa, p_0, s_0\fa)$, that is, a gallery of type the unique minimal gallery from $\sigma = v_0$ to $s_0\fa$, then we would need the label at $\sigma = v_0$ to capture all three of the above distinct types, and their labels.  This does not provide any advantage over the union over $w \in \sW$ described in Theorem \ref{thm:DoubleCosetIntersectionParahoric}.
\end{remark}

We conclude this section by applying the results above to $K = G(\mathcal{O})$, the stabilizer of the origin in $G(F)$.  In the following, for $\lambda \in R^\vee$ we write $\star(\lambda)$ for the star of~$\lambda$ in the standard apartment $\App$, that is, the union of the set of alcoves of $\App$ which contain the vertex $\lambda$.  Recall that for any $\lambda \in R^\vee$, there is a unique alcove $x_\lambda \fa$ in $\star(\lambda)$ which is at minimal distance from the base alcove $\fa$.  Equivalently, $x_\lambda \in \aW$ is the unique minimal length representative of the coset $t^\lambda \sW$ in $\aW/\sW$.

\begin{corollary}\label{cor:DoubleCosetIntersectionK}   Let $\lambda, \mu \in R^\vee$, let $P$ be a standard spherical parabolic subgroup of $G(F)$, and let $y \in \aW$.  Assume $\lambda$ is dominant and let $x_\lambda\fa$ be the alcove in $\star(\lambda)$ closest to $\fa$. Then there is a bijection between the points of the intersection 
	\[Kt^\lambda K \cap (I_P)^y t^\mu K\] 
	and the union over $w\in W_0$ of the set of labeled folded alcove walks of type $\overrightarrow{wx_\lambda}$ from $\fa$ to some alcove in $\star(\mu)$ which are positively folded with respect to the $(P,y)$-chimney. Moreover, 
	\[Kt^\lambda K \cap (I_P)^y t^\mu K \neq \emptyset\]
	if and only if the vertex $\mu$ lies in the shadow $\Sh_{P,y}(\lambda)$. 
\end{corollary}

\noindent Note that, as in Theorems~\ref{thm:DoubleCosetIntersectionParahoric} and~\ref{thm:IntersectionParahoricShadow}, the assumption on $\lambda$ here does not lead to any loss of generality.  The proof below will show why we do not need to specify the final alcove within $\star(\mu)$ of the alcove walks that we consider.

\begin{proof} Since $\lambda$ is dominant, the element $x_\lambda \in \aW$ is $(W_0,W_0)$-reduced, and a minimal gallery from the origin $v_0$ to $\lambda$ of type $\vec{\lambda}$ has first alcove $\fa$ and final alcove $x_\lambda\fa$.  Now $x_\lambda = t^\lambda v$ for some $v \in \sW$, so $K t^\lambda K = K t^\lambda v K = K x_\lambda K$.  Similarly, for any alcove $z\fa$ in $\star(\mu)$, we have $(I_P)^y t^\mu K = (I_P)^y z K$.  We now apply Theorems \ref{thm:DoubleCosetIntersectionParahoric} and \ref{thm:IntersectionParahoricShadow}.
\end{proof}

\section{Examples in rank two}\label{sec:examples}

In this final section, we explore, rather informally, some examples of shadows in affine Coxeter groups of types $\tilde{A}_2$ and $\tilde{C}_2$.  These examples will also illustrate our main recursive result, Theorem~\ref{thm:AlcoveRecursion}.  We did not observe any fundamentally new behavior in type $\tilde{G}_2$ and so omit this case. 

Throughout this section $A$ is a copy of the Coxeter complex of the given affine type, and $R^\vee$ is the corresponding coroot lattice, seen as a subset of the vertices of $A$. Recall from Corollary~\ref{cor:WeylPolytope} that the shadow of a coroot lattice element $\lambda$ is always contained in the $\lambda$-Weyl polytope.  In case the shadow is taken with respect to a chamber at infinity, then by the convexity theorem~\cite{Hitzel}, the shadow of $\lambda$ is equal to the intersection of $R^\vee$ with the $\lambda$-Weyl polytope.  However we will see that for general chimney-induced orientations, shadows may be non-convex, with both ``notches" and ``holes".  We will also observe that some of what is contained in the shadow can be described using ``strings" in root directions.  

We will explain one example in detail in Section~\ref{sec:Avertex}, then sketch several further examples.  

\subsection{Shadows of vertices in type $\tilde{A}_2$} \label{sec:Avertex} 

\begin{figure}[ht]
	\begin{center}
		\resizebox{0.4\textwidth}{!}
		{
			\begin{overpic}{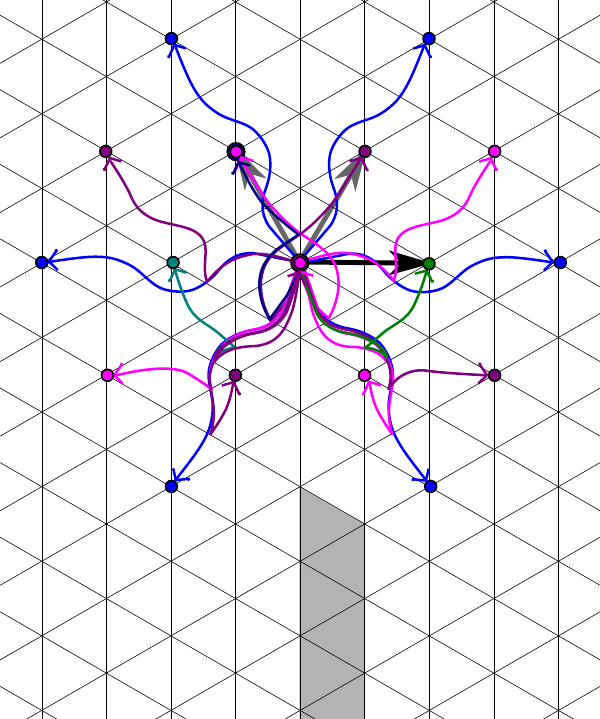}
				\put(52,65){$\alpha$}
				\put(30,73){$\beta$}
				\put(60,96){$\lambda$}
				\put(51,86){\textcolor{blue}{$\gamma_\lambda$}}
			\end{overpic}
		}
		\resizebox{0.4\textwidth}{!}
		{
			\begin{overpic}{AlcoveShadowA2}
				\put(58.5,55.5){$\alpha$}
				\put(48,57){$\beta$}
				\put(28,14){$\x$}
				\put(52,59){$\fa$}
			\end{overpic}
		}
		
		\caption{{Shadows of the coroot lattice element $\lambda$ (on the left) and the alcove $\x$ (on the right), with respect to the chimney represented by the gray shaded region.}}
		\label{fig:chimneyIntroRepeated}
	\end{center}
\end{figure}
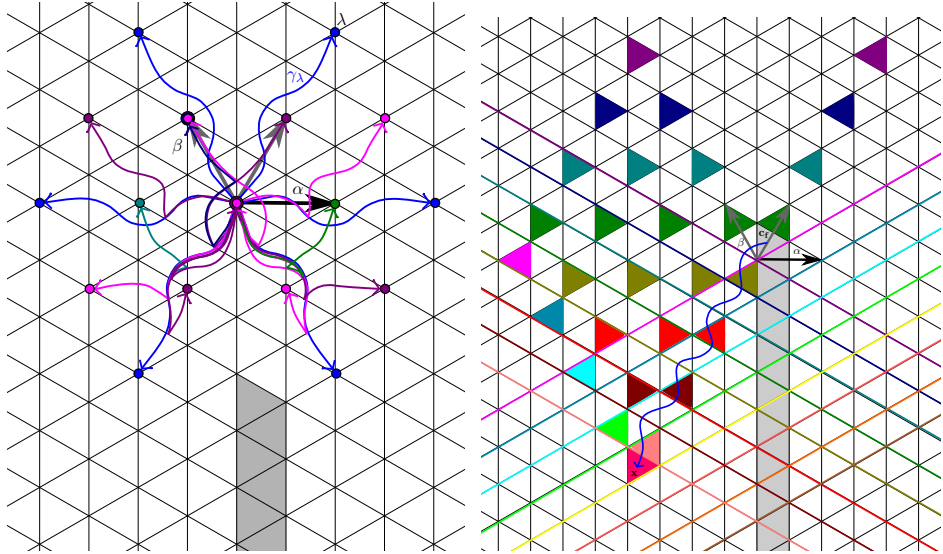

Let us first have a closer look at the picture on the left in Figure~\ref{fig:chimneyIntroRepeated}, which is the same as  Figure~\ref{fig:chimneyIntro} of the introduction.

With $\alpha_1 = \alpha = \alpha^\vee$ and $\alpha_2 = \beta = \beta^\vee$ as given in the figure, let $\lambda=2\alpha+2\beta \in R^\vee$ and define $J=\{1\}$, so that $W_J=\langle s_\alpha\rangle$.  We determine the shadow $\Sh_J(\lambda) = \Sh_{J, \id}(\lambda)$ of $\lambda$ with respect to the $J$-chimney, using the minimal gallery $\gamma_\lambda$ from the origin to $\lambda$.  A $J$-sector representing the $J$-chimney is shown as the gray shaded region in Figure~\ref{fig:chimneyIntroRepeated}.

By definition the elements of $\Sh_J(\lambda)$ are exactly the end-vertices of all positively folded galleries of the same type as $\gamma_\lambda$. We have drawn all those galleries on the left of Figure~\ref{fig:chimneyIntroRepeated}.

We now explain the recursive algorithm given by Theorem~\ref{thm:AlcoveRecursion}.  

\medskip

\noindent {\bf The algorithm:}
For a given coroot lattice element $\lambda$ we obtain from Corollary~\ref{cor:sim-stable-alcove} that  $\Sh_{J}(\lambda)=\Sh_{\y}(\lambda)$ for some alcove $\y$ ``deep enough" inside a $J$-sector. One can then compute $\Sh_{\y}(\lambda)$ as follows: 
\begin{enumerate}
\item Choose a minimal gallery $\gamma_y$ from $\fa$ to $\y$ and denote by $\mathcal{H}_y=(H_1, H_2, \dots ,H_{\ell(y)})$ the sequence of hyperplanes crossed by this gallery (in the order in which they appear). 
\item As the end-vertices of all galleries in $\Gamma_0\define W_0\cdot\gamma_\lambda$ (the blue dots in Figure~\ref{fig:chimneyIntroRepeated}) are in the shadow (each of these is trivially positively folded), put 
\[
\Sh^0_\y(\lambda)\define \sW\cdot\lambda  \subset \Sh_\y(\lambda).
\] 
That is, $\Sh^0_\y(\lambda)$ is the set of end-vertices of galleries in $\Gamma_0$.
\item For $i=1,\dots, \ell(y)$ proceed as follows: 
Let $\Gamma_{i}$ be the set of galleries obtained from $\Gamma_{i-1}$ by introducing (additional) positive folds along the hyperplane $H_i$ in the above sequence $\mathcal{H}_y$, as in the proof of Theorem~\ref{thm:AlcoveRecursion}. Define $\Sh^i_\y(\lambda)$ to be the set of end-vertices of galleries in $\Gamma_i$. Finally we obtain $\Sh_{\y}(\lambda)=\Sh^{\ell(y)}_\y(\lambda)$. 
\end{enumerate}

In the last (repeated) step of the recursive construction, the reflection $r_{H_i}$ will take  end-vertices on the far side of $H_i$ (i.e. vertices that used to be separated from $\fa$ by $H_i$) to end-vertices on the near side (i.e. on the same side of $H_i$ as $\fa$).  Note that a vertex $\mu \in R^\vee$ is on the far side of $H_{\xi,-k}$, with $k$ a non-negative integer, if and only if $\langle \mu,\xi \rangle < -k$.  This procedure (potentially) adds new elements to the shadow.

\medskip
We now describe this recursive procedure in more detail for the case shown on the left of Figure~\ref{fig:chimneyIntroRepeated}. 
Put $s_0\define s_{\alpha + \beta,1}$.  We can choose $\y$ to be of the form $(s_\beta s_\alpha s_0)^N\fa$ for a suitably large $N=N(\lambda)$. We compute the shadow $\Sh_{J}(\lambda)=\Sh_{\y}(\lambda)$ as follows. 

In Step (1), the collection $\mathcal{H}_{y}$ is an initial subsequence of the sequence of hyperplanes
\[
\mathcal{H}_{J} \define (H_{\beta,0}, H_{\alpha+\beta,0}, H_{\beta,-1}, H_{\alpha+\beta,-1}, \dots, H_{\beta,-k}, H_{\alpha+\beta,-k},\dots). 
\]

Step (2) yields  $\Gamma_0:=\sW\cdot\gamma_\lambda$ and
\[
\Sh^0_\y(\lambda)\define \sW\cdot\lambda = \{\pm 2\alpha,\pm 2\beta,\pm 2\alpha+2\beta\} \subset \Sh_{\y}(\lambda).
\]

By the recursive procedure in Step (3) we have to consider the successive hyperplanes in the sequence $\mathcal{H}_J$. The first two such hyperplanes are $H_{\beta, 0}$ and $H_{\alpha+\beta, 0}$.  Applying the corresponding reflections to vertices in $\sW\cdot\lambda$ that are on the far side of these hyperplanes gives us elements once again in $W_0\cdot \lambda$. Therefore adding folds along these two hyperplanes does not add vertices to the shadow. Hence $\Gamma_2=\Gamma_1=\Gamma_0$ and $\Sh^2_\y(\lambda)=\Sh^1_\y(\lambda)=\Sh^0_\y(\lambda)$.  

The next hyperplane is $H_3\define H_{\beta, -1}$. Of the elements of $\sW \cdot \lambda$, only $ 2\alpha$, $-2\beta$, and $-2\alpha-2\beta$ are on the far side of $H_3$. Reflect the tails of the (blue) minimal galleries which have these end-vertices along $H_3$ and add the resulting positively-folded galleries to $\Gamma_2$ to obtain $\Gamma_3$. These new galleries end at $2\alpha + \beta$, $\beta$, and $-2\alpha+\beta$, respectively. Hence $\Sh^3_\y(\lambda) = \Sh^2_\y(\lambda)\cup\{2\alpha + \beta, \beta, -2\alpha+\beta\}$. 

Continuing in this manner we successively consider end-vertices which lie on the far side of the hyperplanes in the sequence $\mathcal{H}_J$ above, and the algorithm terminates when there are no such vertices remaining. 

The shadow consists of all of the coroot lattice points in the $\lambda$-Weyl polytope except for $\pm(\alpha + 2\beta)$ and is hence a strict subset of the intersection of $R^\vee$ with the $\lambda$-Weyl polytope. 
These missing vertices on opposite sides of the $\lambda$-Weyl polytope are an example of what we call ``notches".

\subsection{Shadows of alcoves in type $\tilde{A}_2$} \label{sec:Aalcove} 

We now consider the shadow of an alcove as on the right-hand side of Figure~\ref{fig:chimneyIntroRepeated} (which is the same as Figure~\ref{fig:AlcoveShadowA2} of the introduction).  Let $x \in \aW$ correspond to the (pink) alcove $\x$ at the bottom of the figure and let $\gamma_x:\fa\rightsquigarrow \x$ be the (blue) minimal gallery.  A $J$-sector is shaded in gray. 

The shadow $\Sh_{J}(\x)$ can again be obtained by the algorithm described in Section~\ref{sec:Avertex}, using the same ``deep enough" alcove $\y$ and sequence of successively crossed hyperplanes $\mathcal{H}_J$.  However the sets $\Sh^i_\y(\x)$ now consist of end-alcoves of certain positively folded galleries.  In particular, the ``starting set" $\Sh_\y^0(\x)$ is the set of end-alcoves of all minimal galleries of type~$\gamma_x$ with first alcove $\fa$, and as $\gamma_x$ is the unique such gallery, we have $\Sh_\y^0(\x)=\{ \x \}$.  We have used the same color for a hyperplane $H_i$ in the sequence $\mathcal{H}_J$ and the alcoves in $\Sh_{J}(\x)$ which are (first) obtained by applying the reflection $r_{H_i}$.  

Note that the shadow $\Sh_{J}(\x)$ is not a convex set. More precisely, the shadow of $\x$ is not the intersection of some collection of alcoves with some convex polytope in $A$.  The elements $z \in \aW$ such that $\z$ is in the shadow of $\x$ are, however, always subsets of the Bruhat interval between $x$ and $\id$. See \cite{GraeberSchwer} for the connection between shadows and Bruhat order. 

We also draw attention to the following pattern: for any given root direction, say $\xi$, for which a $\xi$-hyperplane appears in the sequence $\mathcal{H}_J$, there is a first and last index of a $\xi$-hyperplane that appears in $\mathcal{H}_J$ before the algorithm terminates, and all of the indices in between also appear.   For example, $H_{\beta,0}\in \mathcal{H}_J$ and $H_{\beta, -3}\in \mathcal{H}_J$, and all $\beta$-hyperplanes of indices strictly between $0$ and $-3$ are in $\mathcal{H}_J$ as well.  Folding $\gamma_x$ along this collection of hyperplanes yields a ``string" of alcoves in the shadow, with these alcoves differing from each other just by translates $t^\beta$.  There is similarly a ``string" of alcoves in the direction of the root $\alpha + \beta$.  These alcoves can in fact only be obtained by reflecting $\x$ in $(\alpha + \beta)$-hyperplanes, so it is essential to consider non-simple root directions when computing the shadow.

\subsection{Vertex shadows with respect to $J$-chimneys in type $\tilde{C}_2$} \label{sec:CvertexJ} 

\begin{figure}[h]
	\begin{center}
		\resizebox{0.45\textwidth}{!}
		{
			\begin{overpic}{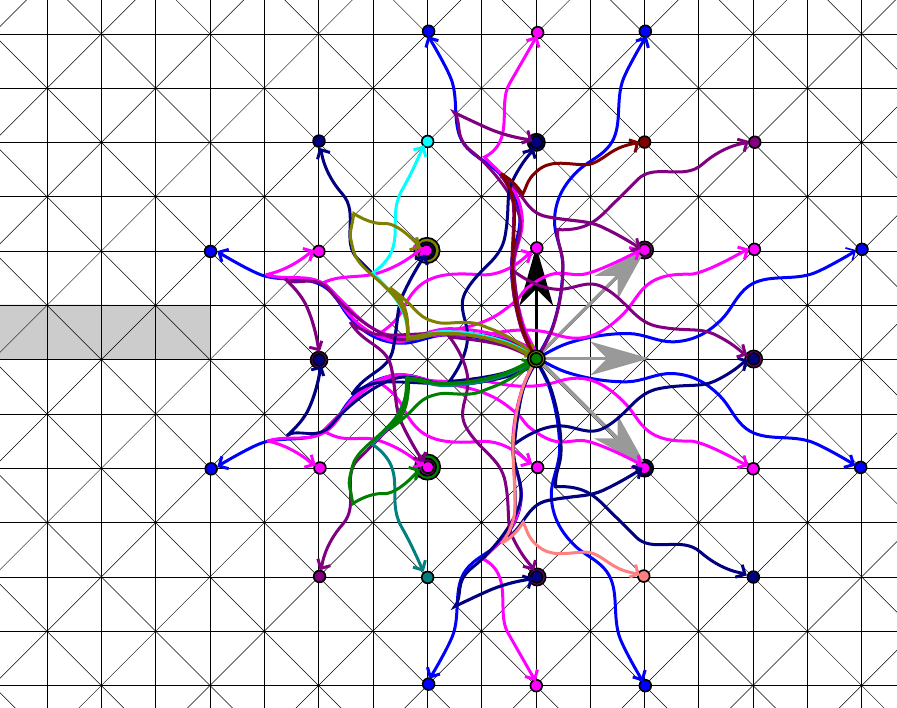}
				\put(60,53){$\alpha$}
				\put(94,53){$\lambda$}
				\put(70,23){$\beta$}
			\end{overpic}
		}
		\resizebox{0.45\textwidth}{!}
		{
			\begin{overpic}{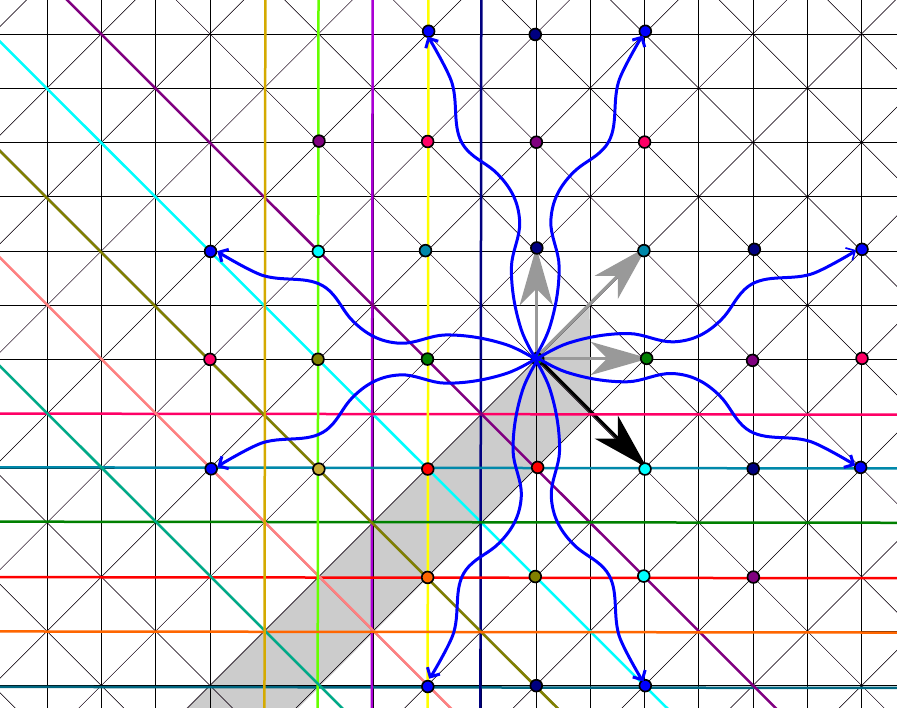}
				\put(60,53){$\alpha$}
				\put(94,53){$\lambda$}
				\put(70,23){$\beta$}
			\end{overpic}
		}
		\caption{Vertex shadows in type $\tilde C_2$ with respect to the chimneys for a short coroot (on the left) and a long coroot (on the right).}
		\label{fig:ShortRootC2}
		\label{fig:LongRootC2}
	\end{center}
\end{figure}
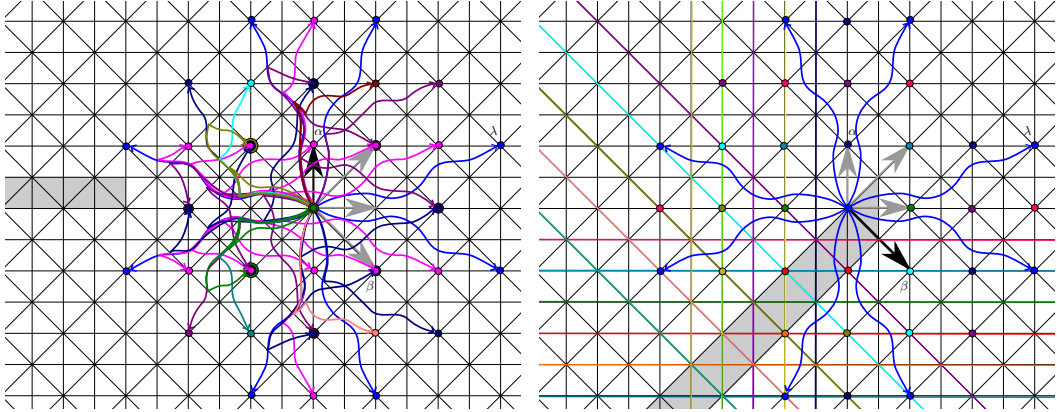

Have a look at the left-hand side of Figure~\ref{fig:ShortRootC2}.
In this example, we depict $\Sh_{J}(\lambda)$ with respect to the $J$-chimney corresponding to the short simple coroot $\alpha^\vee$.  The same recursive procedure as described in Section~\ref{sec:Avertex} can be used to compute this shadow.

We observe that $\Sh_{J,y}(\lambda)$ consists of all the coroot lattice points in the $\lambda$-Weyl polytope {except} for $k(\alpha^\vee + \beta^\vee)$ with $k\in\{\pm 1 \pm 3\}$. That is the set of end-vertices of the galleries depicted. In this case, in addition to ``notches" on opposite sides of the polytope the shadow also has ``holes", consisting of a line of spaced-out coroot lattice points cutting through the polytope. The number of distinct colored-in circles appearing at a coroot lattice point is equal to the multiplicity of this vertex in the shadow.

On the right-hand side of Figure~\ref{fig:LongRootC2}, we depict the shadow of the same $\lambda$ with respect to the $J$-chimney for the long simple coroot $\beta^\vee$.  In this case there are ``notches" but no ``holes", as the shadow consists of all of the coroot lattice points in the $\lambda$-Weyl polytope except for $\pm(4\alpha^\vee + 2\beta^\vee)$.  As on the right of Figure~\ref{fig:chimneyIntroRepeated}, we have not drawn all of the positively folded galleries, but instead illustrate the recursion using color.  The examples in Figure~\ref{fig:ShortRootC2} can be seen to contain ``strings" of shadow elements, in both simple and non-simple root directions.

\subsection{Vertex shadow with respect to a shifted chimney in type $\tilde{C}_2$ }\label{sec:Cshift} 

\begin{figure}[h]
\begin{center}
\resizebox{0.5\textwidth}{!}
{
\begin{overpic}{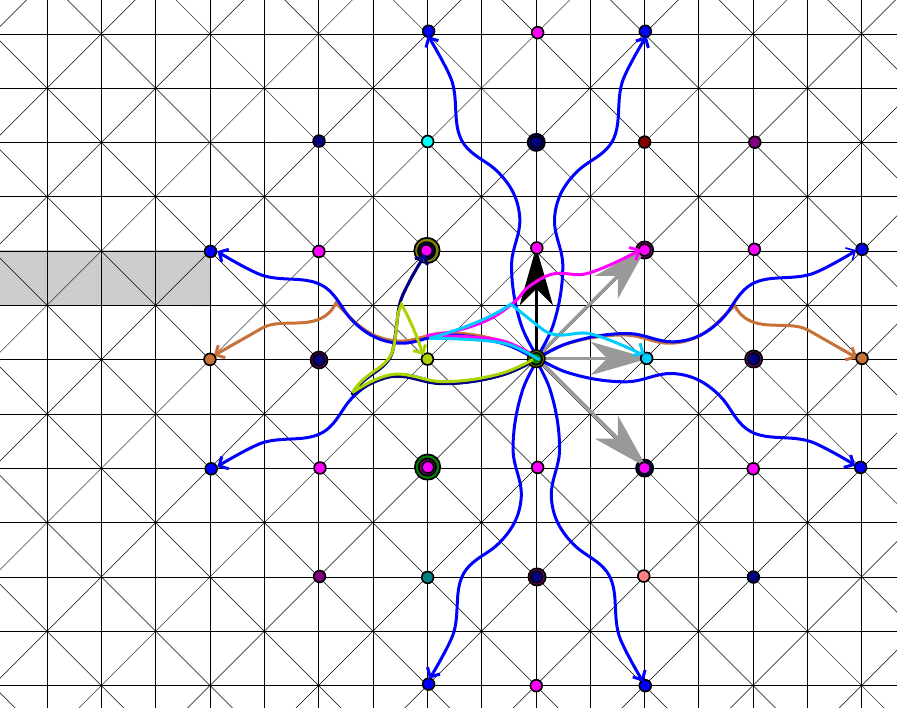}
\put(60,53){$\alpha^\vee$}
\put(70,23){$\beta^\vee$}
\end{overpic}
}
\caption{Vertex shadow with respect to a shifted chimney for a short coroot.}
\label{fig:TranslatedChimneyC2}
\end{center}
\end{figure}
 
 We conclude this examples section by considering the shadow of the same vertex $\lambda = 4\alpha^\vee+3\beta^\vee$ with respect to the chimney that lies between the $\alpha$-hyperplanes of index $1$ and $2$, as shown in Figure~\ref{fig:TranslatedChimneyC2}.

To find the shadow in this case, we may initially cross the same sequence of hyperplanes as for the $J$-chimney lying between $H_{\alpha,0}$ and $H_{\alpha,1}$, as on the left of Figure~\ref{fig:ShortRootC2}.  Thus we initially obtain the same collection of vertices in the shadow as before.  We then cross $H_{\alpha,1}$ to reach an alcove which lies deep in the shaded region depicted in Figure~\ref{fig:TranslatedChimneyC2}.  At this last step, we obtain new folded galleries, by reflecting in $H_{\alpha,1}$. These additional galleries (and their predecessors) are drawn in Figure~\ref{fig:TranslatedChimneyC2}.   Therefore the shadow in this case contains more vertices than before.  And in fact, the shadow now equals the intersection of $R^\vee$ with the $\lambda$-Weyl polytope (shown with multiplicity).

\section{Conflict of interest}

The authors declare that they have no conflict of interest.

\renewcommand{\refname}{Bibliography}
\bibliography{bibliography}
\bibliographystyle{alpha}

\end{document}